\catcode`\@=11

\font\tenmsa=msam10 \font\sevenmsa=msam7 \font\fivemsa=msam5
\font\tenmsb=msbm10
\font\sevenmsb=msbm7 \font\fivemsb=msbm5 \newfam\msafam \newfam\msbfam
\textfont\msafam=\tenmsa \scriptfont\msafam=\sevenmsa
\scriptscriptfont\msafam=\fivemsa \textfont\msbfam=\tenmsb
\scriptfont\msbfam=\sevenmsb \scriptscriptfont\msbfam=\fivemsb

\def\hexnumber@#1{\ifnum#1<10 \number#1\else \ifnum#1=10 A\else\ifnum#1=11
 B\else\ifnum#1=12 C\else \ifnum#1=13 D\else\ifnum#1=14 E\else\ifnum#1=15
 F\fi\fi\fi\fi\fi\fi\fi}

\def\msa@{\hexnumber@\msafam} \def\msb@{\hexnumber@\msbfam}
\mathchardef\boxdot="2\msa@00 \mathchardef\boxplus="2\msa@01
\mathchardef\boxtimes="2\msa@02 \mathchardef\square="0\msa@03
\mathchardef\blacksquare="0\msa@04 \mathchardef\centerdot="2\msa@05
\mathchardef\lozenge="0\msa@06 \mathchardef\blacklozenge="0\msa@07
\mathchardef\circlearrowright="3\msa@08 \mathchardef\circlearrowleft="3\msa@09
\mathchardef\rightleftharpoons="3\msa@0A
\mathchardef\leftrightharpoons="3\msa@0B \mathchardef\boxminus="2\msa@0C
\mathchardef\Vdash="3\msa@0D \mathchardef\Vvdash="3\msa@0E
\mathchardef\vDash="3\msa@0F \mathchardef\twoheadrightarrow="3\msa@10
\mathchardef\twoheadleftarrow="3\msa@11 \mathchardef\leftleftarrows="3\msa@12
\mathchardef\rightrightarrows="3\msa@13 \mathchardef\upuparrows="3\msa@14
\mathchardef\downdownarrows="3\msa@15 \mathchardef\upharpoonright="3\msa@16
 \mathchardef\downharpoonright="3\msa@17
\mathchardef\upharpoonleft="3\msa@18 \mathchardef\downharpoonleft="3\msa@19
\mathchardef\rightarrowtail="3\msa@1A \mathchardef\leftarrowtail="3\msa@1B
\mathchardef\leftrightarrows="3\msa@1C \mathchardef\rightleftarrows="3\msa@1D
\mathchardef\Lsh="3\msa@1E \mathchardef\Rsh="3\msa@1F
\mathchardef\rightsquigarrow="3\msa@20
\mathchardef\leftrightsquigarrow="3\msa@21 \mathchardef\looparrowleft="3\msa@22
\mathchardef\looparrowright="3\msa@23 \mathchardef\circeq="3\msa@24
\mathchardef\succsim="3\msa@25 \mathchardef\gtrsim="3\msa@26
\mathchardef\gtrapprox="3\msa@27 \mathchardef\multimap="3\msa@28
\mathchardef\therefore="3\msa@29 \mathchardef\because="3\msa@2A
\mathchardef\doteqdot="3\msa@2B 
\mathchardef\traceiangleq="3\msa@2C \mathchardef\precsim="3\msa@2D
\mathchardef\lesssim="3\msa@2E \mathchardef\lessapprox="3\msa@2F
\mathchardef\eqslantless="3\msa@30 \mathchardef\eqslantgtr="3\msa@31
\mathchardef\curlyeqprec="3\msa@32 \mathchardef\curlyeqsucc="3\msa@33
\mathchardef\preccurlyeq="3\msa@34 \mathchardef\leqq="3\msa@35
\mathchardef\leqslant="3\msa@36 \mathchardef\lessgtr="3\msa@37
\mathchardef\backprime="0\msa@38 \mathchardef\risingdotseq="3\msa@3A
\mathchardef\fallingdotseq="3\msa@3B \mathchardef\succcurlyeq="3\msa@3C
\mathchardef\geqq="3\msa@3D \mathchardef\geqslant="3\msa@3E
\mathchardef\gtrless="3\msa@3F \mathchardef\sqsubset="3\msa@40
\mathchardef\sqsupset="3\msa@41
\mathchardef\trianglelefteq="3\msa@45 \mathchardef\bigstar="0\msa@46
\mathchardef\between="3\msa@47 \mathchardef\blacktriangledown="0\msa@48
\mathchardef\blacktriangleright="3\msa@49
\mathchardef\blacktriangleleft="3\msa@4A
\mathchardef\blacktriangle="0\msa@4E \mathchardef\triangledown="0\msa@4F
\mathchardef\eqcirc="3\msa@50 \mathchardef\lesseqgtr="3\msa@51
\mathchardef\gtreqless="3\msa@52 \mathchardef\lesseqqgtr="3\msa@53
\mathchardef\gtreqqless="3\msa@54 \mathchardef\Rrightarrow="3\msa@56
\mathchardef\Lleftarrow="3\msa@57 \mathchardef\veebar="2\msa@59
\mathchardef\barwedge="2\msa@5A \mathchardef\doublebarwedge="2\msa@5B
\mathchardef\angle="0\msa@5C \mathchardef\measuredangle="0\msa@5D
\mathchardef\sphericalangle="0\msa@5E \mathchardef\varpropto="3\msa@5F
\mathchardef\smallsmile="3\msa@60 \mathchardef\smallfrown="3\msa@61
\mathchardef\Subset="3\msa@62 \mathchardef\Supset="3\msa@63
\mathchardef\Cup="2\msa@64  \mathchardef\Cap="2\msa@65
 \mathchardef\curlywedge="2\msa@66
\mathchardef\curlyvee="2\msa@67 \mathchardef\leftthreetimes="2\msa@68
\mathchardef\rightthreetimes="2\msa@69 \mathchardef\subseteqq="3\msa@6A
\mathchardef\supseteqq="3\msa@6B \mathchardef\bumpeq="3\msa@6C
\mathchardef\Bumpeq="3\msa@6D \mathchardef\lll="3\msa@6E 
\mathchardef\ggg="3\msa@6F  \mathchardef\circledS="0\msa@73
\mathchardef\pitchfork="3\msa@74 \mathchardef\dotplus="2\msa@75
\mathchardef\backsim="3\msa@76 \mathchardef\backsimeq="3\msa@77
\mathchardef\complement="0\msa@7B \mathchardef\intercal="2\msa@7C
\mathchardef\circledcirc="2\msa@7D \mathchardef\circledast="2\msa@7E
\mathchardef\circleddash="2\msa@7F \def\ulcorner{\delimiter"4\msa@70\msa@70 }
\def\urcorner{\delimiter"5\msa@71\msa@71 }
\def\llcorner{\delimiter"4\msa@78\msa@78 }
\def\lrcorner{\delimiter"5\msa@79\msa@79 } \def\yen{\mathhexbox\msa@55 }
\def\checkmark{\mathhexbox\msa@58 } \def\circledR{\mathhexbox\msa@72 }
\def\maltese{\mathhexbox\msa@7A } \mathchardef\lvertneqq="3\msb@00
\mathchardef\gvertneqq="3\msb@01 \mathchardef\nleq="3\msb@02
\mathchardef\ngeq="3\msb@03 \mathchardef\nless="3\msb@04
\mathchardef\ngtr="3\msb@05 \mathchardef\nprec="3\msb@06
\mathchardef\nsucc="3\msb@07 \mathchardef\lneqq="3\msb@08
\mathchardef\gneqq="3\msb@09 \mathchardef\nleqslant="3\msb@0A
\mathchardef\ngeqslant="3\msb@0B \mathchardef\lneq="3\msb@0C
\mathchardef\gneq="3\msb@0D \mathchardef\npreceq="3\msb@0E
\mathchardef\nsucceq="3\msb@0F \mathchardef\precnsim="3\msb@10
\mathchardef\succnsim="3\msb@11 \mathchardef\lnsim="3\msb@12
\mathchardef\gnsim="3\msb@13 \mathchardef\nleqq="3\msb@14
\mathchardef\ngeqq="3\msb@15 \mathchardef\precneqq="3\msb@16
\mathchardef\succneqq="3\msb@17 \mathchardef\precnapprox="3\msb@18
\mathchardef\succnapprox="3\msb@19 \mathchardef\lnapprox="3\msb@1A
\mathchardef\gnapprox="3\msb@1B \mathchardef\nsim="3\msb@1C
\mathchardef\napprox="3\msb@1D
\mathchardef\nsupseteqq="3\msb@23 \mathchardef\subsetneqq="3\msb@24
\mathchardef\supsetneqq="3\msb@25
\mathchardef\supsetneq="3\msb@29 \mathchardef\nsubseteq="3\msb@2A
\mathchardef\nsupseteq="3\msb@2B \mathchardef\nparallel="3\msb@2C
\mathchardef\nmid="3\msb@2D \mathchardef\nshortmid="3\msb@2E
\mathchardef\nshortparallel="3\msb@2F \mathchardef\nvdash="3\msb@30
\mathchardef\nVdash="3\msb@31 \mathchardef\nvDash="3\msb@32
\mathchardef\nVDash="3\msb@33 \mathchardef\ntrianglerighteq="3\msb@34
\mathchardef\ntrianglelefteq="3\msb@35 \mathchardef\ntriangleleft="3\msb@36
\mathchardef\ntriangleright="3\msb@37 \mathchardef\nleftarrow="3\msb@38
\mathchardef\nrightarrow="3\msb@39 \mathchardef\nLeftarrow="3\msb@3A
\mathchardef\nRightarrow="3\msb@3B \mathchardef\nLeftrightarrow="3\msb@3C
\mathchardef\nleftrightarrow="3\msb@3D \mathchardef\divideontimes="2\msb@3E
\mathchardef\varnothing="0\msb@3F \mathchardef\nexists="0\msb@40
\mathchardef\mho="0\msb@66 \mathchardef\thorn="0\msb@67
\mathchardef\beth="0\msb@69 \mathchardef\gimel="0\msb@6A
\mathchardef\daleth="0\msb@6B \mathchardef\lessdot="3\msb@6C
\mathchardef\gtrdot="3\msb@6D \mathchardef\ltimes="2\msb@6E
\mathchardef\rtimes="2\msb@6F \mathchardef\shortmid="3\msb@70
\mathchardef\shortparallel="3\msb@71 \mathchardef\smallsetminus="2\msb@72
\mathchardef\thicksim="3\msb@73 \mathchardef\thickapprox="3\msb@74
\mathchardef\approxeq="3\msb@75 \mathchardef\succapprox="3\msb@76
\mathchardef\precapprox="3\msb@77 \mathchardef\curvearrowleft="3\msb@78
\mathchardef\curvearrowright="3\msb@79 \mathchardef\digamma="0\msb@7A
\mathchardef\varkappa="0\msb@7B \mathchardef\hslash="0\msb@7D
\mathchardef\hbar="0\msb@7E \mathchardef\backepsilon="3\msb@7F
\def\Bbb{\ifmmode\let\next\Bbb@\else
\def\next{\errmessage{Use \string\Bbb\space only in math mode}}\fi\next}
\def\Bbb@#1{{\Bbb@@{#1}}} \def\Bbb@@#1{\fam\msbfam#1}

\catcode`\@=\active

\def\del{\partial}
\def\CL{\hbox{{$\cal L$}}}

\def\CR{\hbox{{$\cal R$}}}

 \def\CZ{\hbox{{$\cal Z$}}}

\def\cu{{\upsilon}} % used for special element u
\def\cv{{\vartheta}} % used for special element v

\def\cg{\hbox{{\sl g}}} % used for Lie algebra 'gothic g'
\def\lform{\hbox{$\sqcup$}\llap{\hbox{$\sqcap$}}}
\def\h{{{1\over2}}}

\def\R{{\Bbb R}}
\def\C{{\Bbb C}}
\def\Z{{\Bbb Z}}

\def\eps{{\epsilon}}

\def\dcross{{\bowtie}}
\def\codcross{{\blacktriangleright\!\!\blacktriangleleft}}
\def\rbiprod{{\cdot\kern-.33em\triangleright\!\!\!<}}
\def\lbiprod{{>\!\!\!\triangleleft\kern-.33em\cdot\, }}

\def\tens{\mathop{\otimes}}

\def\la{{\triangleright}}

\def\isom{{\cong}}

\def\ev{{\rm ev}}
\def\coev{{\rm coev}}

\def\id{{\rm id}}

\def\<{\langle}
\def\>{\rangle}

\def\equad{\kern -1.7em}

\def\eqn#1#2{\begin{equation}#2\label{#1}\end{equation}}

\def\haj#1{{\mathaccent20 {#1}}}
\def\Vhaj{{V\haj{\ }}}

\def\o{{}_{\scriptscriptstyle(1)}}
\def\t{{}_{\scriptscriptstyle(2)}}
\def\th{{}_{\scriptscriptstyle(3)}}
\def\fo{{}_{\scriptscriptstyle(4)}}

\def\bo{{}^{\bar{\scriptscriptstyle(1)}}}
\def\bt{{}^{\bar{\scriptscriptstyle(2)}}}

\def\und#1{{\underline {#1}}}

\def\uo{{{}^{\scriptscriptstyle(1)}}}
\def\ut{{{}^{\scriptscriptstyle(2)}}}

\def\umo{{{}^{\scriptscriptstyle-(1)}}}
\def\umt{{{}^{\scriptscriptstyle-(2)}}}

\def\Bo{{{}_{\und{\scriptscriptstyle(1)}}}}
\def\Bt{{{}_{\und{\scriptscriptstyle(2)}}}}

\def\text#1{\mbox{\rm #1}}
\def\note#1{}

\def\blacksquare{{\lform}}%AMS Tex Fakes
\def\frac#1#2{{{#1\over#2}}}

\def\proof{\goodbreak\noindent{\bf Proof\quad}}

\def\endproof{{\ $\lform$}\bigskip }

\def\align#1{\begin{eqnarray*}#1\end{eqnarray*}}
\def\cmath#1{\[\begin{array}{c} #1 \end{array}\]}
\def\ceqn#1#2{\begin{equation}\label{#1}\begin{array}{c}#2
\end{array}\end{equation}}

\def\vect{{\bf t}}\def\vecs{{\bf s}}
\def\vecu{{\bf u}}\def\vecx{{\bf x}}\def\vecp{{\bf p}}
\def\vecl{{\bf l}}

\def\vecm{{\bf m}}

\def\bramat{R_{21}\vecu_1 R\vecu_2=\vecu_2 R_{21}\vecu_1 R}
\def\eucmat{R_{21}\vecx_1 \vecx_2=\vecx_2\vecx_1 R}

\def\baro{{{}_{\bar{(1)}}}}
\def\bart{{{}_{\bar{(2)}}}}

\documentstyle[11pt,epsf]{article} \textheight 23.6cm \textwidth 16cm
\evensidemargin
0in \topskip 28pt

\newtheorem{lemma}{Lemma}[section] \newtheorem{propos}[lemma]{Proposition}
 \newtheorem{theorem}[lemma]{Theorem}
\newtheorem{cor}[lemma]{Corollary} \newtheorem{corol}[lemma]{Corollary}
\newtheorem{defin}[lemma]{Definition}

\begin{document}\baselineskip 20pt

{\ }\hskip 4.7in DAMTP/95-11
\vspace{.2in}

\begin{center} {\LARGE  QUASI-$*$ STRUCTURE ON $q$-POINCARE ALGEBRAS} \\
\baselineskip 13pt
{\ } {\ }\\ S.  Majid\footnote{Royal Society University Research Fellow and
Fellow of Pembroke College, Cambridge}\\{\ }\\ Department of Applied
Mathematics \& Theoretical Physics\\ University of Cambridge, Cambridge CB3
9EW, UK\\
+\\
Department of Mathematics, Harvard University\\
1 Oxford Street, Cambridge, MA02138, USA\footnote{During the calendar years
1995+1996}
\end{center}

\vspace{10pt} \begin{quote}\baselineskip 12pt
\noindent{\bf Abstract} We use braided groups to introduce a theory of
$*$-structures on general inhomogeneous quantum groups, which we formulate as
{\em quasi-$*$} Hopf algebras. This allows the construction of the tensor
product of unitary representations up to a quantum cocycle isomorphism, which
is a novel feature of the inhomogeneous case. Examples include $q$-Poincar\'e
quantum group enveloping algebras in $R$-matrix form appropriate to the
previous $q$-Euclidean and $q$-Minkowski spacetime algebras $\eucmat$ and
$\bramat$. We obtain unitarity of the  fundamental differential
representations. We show further that the Euclidean and Minkowski Poincar\'e
quantum groups are twisting equivalent by a another quantum cocycle.

\end{quote}
\vspace{8pt}
\centerline{March 1995}

\baselineskip 15pt

\tableofcontents

\baselineskip 22pt

\section{Introduction}

As well as specific roles in physical systems, quantum groups in recent years
have motivated a quite general and systematic development of a kind of
$q$-deformed geometry. The basic algebraic ingredients are quite well
understood by now, at least for the geometry of $q$-deformed compact groups
(typically quantum groups\cite{Dri}\cite{Jim:dif}) and $q$-deformed linear
spaces (typically the more novel braided groups introduced by the
author\cite{Ma:exa}). Not understood, however, is the full story regarding the
role of the $*$-structure or complex conjugation in this $q$-geometry. This is
obviously important for contact with physics, where $q$-deformed field theories
are expected to be of interest either as providing a regularisation of
infinities as poles at $q=1$\cite{Ma:reg} or as effective theories modelling
the feedback of quantum or other effects on geometry at the Planck
scale\cite{Ma:pla}\cite{Kem:fie}. An understanding of the $*$-structure is
needed also for better contact with other $C^*$-algebra approaches to
non-commutative geometry\cite{Con:non}.

The present work follows on from a previous one\cite{Ma:star} where we studied
the $*$-structures on linear braided groups. Now we combine the considerations
there with a previous general construction\cite{Ma:poi} for inhomogeneous
quantum groups in order to develop a theory of $*$-structures on these.
Many authors have considered $q$-Poincar\'e and  other inhomogeneous quantum
groups as a key step for $q$-deformed physics but found that they do not (in
the examples that concern us here) have $*$-structures obeying the usual
axioms\cite{Wor:com} of a Hopf $*$-algebra. In physical terms it means that the
`unitary' ($*$-preserving) representations of these inhomogeneous quantum
groups are not closed under tensor product. This is problematic and tells us
that we need a more radical $q$-deformation of the concept of `unitarity' as
well. We propose in the present paper a solution to this long-standing problem.
It was announced briefly at the end of \cite{Ma:clau94} and is developed now in
detail. In physical terms, we will see that quantum deformation introduces a
kind of `anomaly' in the sense of a cocycle governing the breakdown of
unitarity. Actually, something like this is to be expected because if we view
$q$-deformation as a regularisation scheme for field theory, we do have to
recover anomalies in the cases where they exist. For example, unlike
dimensional regularisation there is no problem with the $\eps$ tensor and we
have to expect a problem in different quarter to generate the $U(1)$ axial
anomaly.

Namely, we introduce the notion of a {\em quasi}-$*$ Hopf algebra, which is a
Hopf algebra  where the algebra part is a $*$-algebra and the coalgebra
$\Delta,\eps$ obeys
\cmath{ (*\tens *)\circ\Delta\circ *=\CR^{-1}(\tau\circ\Delta\ )\CR,\quad
\overline{\eps(\ )}=\eps\circ * ,\quad \CR^{*\tens *}=\CR_{21} \\
(\id\tens\Delta)\CR=\CR_{13}\CR_{12},\quad
(\Delta\tens\id)\CR=\CR_{13}\CR_{23}}
for some element $\CR$ in the tensor square. We use the usual notations as in
\cite{Swe:hop}\cite{Ma:qua} or the forthcoming text\cite{Ma:book}. The above
axioms are a generalisation both of the usual axioms of a Hopf $*$-algebra and
of Drinfeld's axioms of a quasitriangular Hopf algebra. They reduce to one iff
they reduce to the other, in which case they reduce to a quasitriangular Hopf
$*$-algebra of real type. This formulation covers the examples that interest us
when $q$ is real: one can also consider a different framework suitable for
$\CR^{*\tens *}=\CR^{-1}$ though we do not do so here explicitly. We arrive at
these axioms in Section~2 where we study in detail the abstract construction of
inhomogeneous quantum groups from the `braided geometry' approach based on a
process of {\em bosonisation} introduced by the author in \cite{Ma:bos}. This
associated to any braided group (in our case the linear or `momentum' sector of
the inhomogeneous quantum group) an ordinary quantum group by a certain
semidirect product construction. We also study in Section~2 how the
bosonisation changes under twisting by a quantum cocycle, cf the ideas of
Drinfeld in \cite{Dri:qua}.

In the `braided approach' to $q$-deformed geometry we begin with braided group
deformations of $\R^n$ as the basic objects. The braiding (and $q$) enters in a
way that is conceptually different from other approaches to non-commutative
geometry, namely as {\em braid-statistics} with which we explicitly endow the
co-ordinates of the $q$-deformed $\R^n$. This is the key difference between the
more familiar quantum groups (which are bosonic objects) and the new braided
groups. For a full introduction to the latter we refer to
\cite{Ma:introp}\cite{Ma:introm}\cite{Ma:varen} as well as the original
works\cite{Ma:exa}\cite{Ma:skl} etc. There is a solid theory of braided
matrices, braided linear algebra, braided addition, braided differentiation,
`Poincar\'e' quantum groups, differential forms, epsilon tensors and
integration developed in a series of
papers\cite{Ma:exa}\cite{Ma:lin}\cite{Ma:poi}\cite{Ma:fre}\cite{Ma:eps}
\cite{KemMa:alg}. There are also natural candidates for $q$-Minkowski
\cite{Ma:mec}\cite{Mey:new}\cite{MaMey:bra} and $q$-Euclidean\cite{Ma:euc}
spacetime algebras within this programme, making contact too with earlier
pioneering work in \cite{CWSSW:lor}\cite{CWSSW:ten}\cite{OSWZ:def} where
the same spacetime algebras were proposed directly by other means. The
braided approach extended the latter works and put them, moreover, into a
general R-matrix form as part of uniform theory of braided spaces. The
present work on quasi-$*$ Hopf algebras applies to all the inhomogeneous
quantum groups in this approach\cite{Ma:poi} for which suitable reality
properties are met. This includes the $q$-Euclidean group of motions
$\R_q^n\lbiprod U_q(so_n)$ in any dimension, and also to the 4-dimensional
Minkowski version. These and other examples are described in detail in
Section~3. We also show that the inhomogeneous quantum groups for the
4-dimensional Euclidean and Minkowski cases are related by twisting,
extending the `quantum wick rotation' in \cite{Ma:euc}.

In Section~4 we look at the representation theory of quasi-$*$ Hopf algebras
and inhomogeneous quantum groups from the point of view of our braided
approach. This approach brings out the deeper meaning of the role of $*$ in
braided geometry as a combined conjugation-braid reversal symmetry of our
constructions. The seeds of this idea are already implicit in the earliest
works on quantum groups, where the conjugate $\vecl^\pm$ generators of a
quantum enveloping algebra\cite{FRT:lie} are associated to a universal R-matrix
and its inverse-transpose respectively.  In braided geometry all our
constructions are done by means of braid and tangle diagrams representing the
`flow' of algebraic information, with the braiding $\Psi=\epsfbox{braid.eps}$
playing the role of usual transpositions between independent symbols in
ordinary mathematics, or of super-transposition between independent symbols in
super-constructions. This means, however, that in making braided versions of
classical constructions we have to make choices of under- or over- braid
crossing. Whatever construction we do, we could make a parallel one  in the
braided category with `conjugate' or inverse braiding where the role of under-
and over- braid crossings is reversed\cite{Ma:introp}.
The idea is that in the braided approach our usual classical geometry {\em
splits}  into {\em two} braided versions related by a combined braid-reversal
and $*$ symmetry. This is a new phenomenon not visible classically or even in
super geometry, where $\Psi^2=\id$. Note that these braids do not live in
physical space but in the 3-dimensional `lexicographical space' in which we
write our mathematical constructions diagrammatically.

The idea of $*$ mapping between two versions of $q$-deformed geometry rather
than being (as more usual) a property of one system, is evident in the theory
of $*$-structure on linear spaces developed in \cite{Ma:star}, where, for
example, we viewed $*:\Omega_L\to \Omega_R$ between left and right versions of
the $q$-deformed exterior algebra. We have the same phenomenon for $q$-deformed
Poincar\'e quantum group function algebras. In the present context of the
Poincar\'e enveloping algebra quantum groups it means that we have two natural
coproducts, connected by $*$. This means in turn that we have two natural
representations of each $q$-Poincar\'e enveloping algebra by braided
differentiation on the $q$-spacetime co-ordinates. We can use either the
differentials $\del$ from \cite{Ma:fre} or `conjugate' ones $\bar \del$ defined
with the inverse braiding in their braided-Leibniz rules. This makes contact
with examples of `conjugate' derivatives in \cite{OgiZum:rea} and elsewhere, as
well as with \cite{Fio:man} where a `zero-dimensional' toy model for a
Poincar\'e quantum group with two coproducts related by $*$ was considered. The
new feature in our approach is that these derivatives are constructed quite
systematically from the braided coaddition (from the left and right) and as
such we now obtain a precise and completely general understanding of how they
are related by $*$ and with each other. The first main result, in Section~2, is
that these two conjugate representations are isomorphic, being intertwined as
\[ \und S\circ\del^i=-\bar\del^i\circ\und S\]
by the braided-antipode $\und S$ or {\em quantum parity operator} $\vecx\to
-\vecx$ on our $q$-spacetime co-ordinates. The second main result, in
Section~4, is a general construction for a sesquilinear form or `inner product'
with respect to which the fundamental representation by translation and
rotation is unitary. A general  feature is that it is no longer exactly
conjugate-symmetric but only, in case of $q$-Euclidean and $q$-Minkowski
spaces, up to a power of $q$.  We obtain in principle a braided version of the
$L^2$ inner product on braided linear spaces, such that $\del$ and $\bar\del$
are mutually adjoint (or such that $\del$ is self-adjoint up to $\und S$). The
computation of such inner products and development of the attendant `braided
analysis' are a direction for further work.

We note that some of the R-matrix formulae  in the present paper can (once
found) be partially verified by direct calculations using the quantum
Yang-Baxter equations (QYBE) many times. This would not, however, check the
various other non-R-matrix relations among quantum group generators $\vect$,
$\vecl^\pm$  etc. in the notation of \cite{FRT:lie}. For a rigorous treatment
that includes all such details automatically, the abstract setting which we use
in the present paper is really needed. We use it in Sections~2 and~4 where we
develop elements of the abstract bosonisation theory from \cite{Ma:bos}, which
in turn ensures consistency with respect to all such additional relations in
the R-matrix formulae presented in Sections~3 and~5.

\subsection*{Preliminaries}

We assume that the reader is familiar with the definition of a quasitriangular
Hopf algebra $H$ in \cite{Dri} with quasitriangular structure or `universal
R-matrix' $\CR=\CR\uo\tens\CR\ut$ in $H\tens H$ (summation of terms implicit),
and the
dual notion of a dual-quasitriangular Hopf algebra $A$ with
dual-quasitriangular structure or `universal R-matrix functional' $\CR:A\tens
A\to \C$ in \cite{Ma:eul}\cite{Ma:bg}. The latter is characterised by the
axioms
\ceqn{dqua}{ \CR(a\tens bc)=\CR(a\o\tens c)\CR(a\t\tens b),\quad
\CR(ab\tens c)=\CR(a\tens c\o)\CR(b\tens c\t)\\ b\o a\o \CR(a\t\tens
b\t)=\CR(a\o\tens b\o)a\t b\t.}
The coproducts are denoted $\Delta a=a\o\tens a\t$, etc. (summation implicit).

It is well-known that the category of representations (modules) of a
quasitriangular Hopf algebra form a braided category with braiding
$\Psi=\epsfbox{braid.eps}$ given by acting via $\CR$ and then making the usual
transposition of the underlying vector
spaces\cite{Dri}\cite{Ma:qua}\cite{ResTur:rib}. By braided category we mean a
collection of objects with a tensor product which is associative and
commutative up to isomorphism. We suppress the associativity isomorphism (which
is trivial in our examples) and write the commutativity isomorphism as the
braiding $\Psi$. There are various coherence axioms between these structures to
the effect that the rules for working in a braided category are the obvious
ones suggested by the braid-crossing notation. The category of
corepresentations (comodules) of a dual-quasitriangular Hopf algebra likewise
form a braided category, with $\Psi$ given by coacting on each comodule,
evaluating the relevant outputs against $\CR$ and making the usual
transposition of the underlying vector spaces vector
spaces\cite{Ma:eul}\cite{Ma:bg}. Introductions are in \cite{Ma:introm}.

The notion of duality of Hopf algebras is best handled for our purposes as
duality pairing\cite{Ma:qua} between two Hopf algebras rather than regarding
one as a subspace of linear functionals on the other. The
product of one maps as usual to the coproduct of the other,
\eqn{hopfdual}{\<h,ab\>=\<h\o,a\>\<h\t,b\>,\quad
\<hg,a\>=\<h,a\o\>\<g,a\t\>,\quad \<h,Sa\>=\<Sh,a\>}
for all $a,b\in A$, $h,g\in H$. The axioms of a Hopf $*$-algebra have been
studied extensively by Woronowicz\cite{Wor:com} and are that our Hopf algebra
is a $*$-algebra and
\eqn{hopf*}{(*\tens *)\circ\Delta=\Delta\circ*,\quad \overline{\eps(\
)}=\eps\circ*,\quad *\circ S=S^{-1}\circ *.}
When we have a quasitriangular structure, it is natural to require $\CR^{*\tens
*}=\CR_{21}$ or $\CR^{*\tens *}=\CR^{-1}$ as explained in
\cite{Dri:alm}\cite{Ma:mec} among other places. The first type is called {\em
real quasitriangular}.

We also assume that the reader is familiar with the basic notion of a braided
group $B$ or braided-Hopf algebra\cite{Ma:bra}\cite{Ma:bg}. Introductions are
in \cite{Ma:introp}\cite{Ma:introm}\cite{Ma:varen}. A braided group $B$ is
like a quantum group but the coproduct $\und\Delta: B\to B\und\tens B$ maps to
the {\em braided tensor product} where the two factors in the tensor product do
not commute. Instead they enjoy mutual {\em braid statistics}. In mathematical
terms we have a Hopf algebra in a braided category, where the braiding $\Psi$
between any two objects determines their mutual braided tensor product algebra.
There are also more usual axioms for a braided antipode $\und S:B\to B$ and
counit $\und\eps:B\to\C$. The diagrammatic way of working with braided groups
consists of writing all maps as arrows generally pointing downwards. We write
tensor products of objects by horizontal juxtaposition,
$\Psi=\epsfbox{braid.eps}$ as a braid crossing and $\Psi^{-1}$ as the reversed
braid crossing. We write other morphisms as nodes with the appropriate valency.
So the product morphism of a braided group is $\cdot=\epsfbox{prodfrag.eps}$
and the coproduct morphism is $\und\Delta=\epsfbox{deltafrag.eps}$.
Functoriality of the braiding is expressed as being allowed to pull nodes
through braid crossings
as if they are beads on the string or tangle. The appropriate coherence theorem
ensures that these rules are consistent and that `topologically equivalent'
diagrams correspond to the same algebraic operations. This kind of `braided
algebra' appeared in \cite{Ma:bos} and is a characteristic feature of the
theory of braided groups\cite{Ma:introp}\cite{Ma:introm}.

In practice, one can also work with braided groups without being too
categorical: one can simply specify every time the required cross-relations or
braid-statistics in the braided tensor product $B\und\tens B$. We write
$b\equiv b\tens 1$ for elements in the first factor and $b'\equiv1\tens 1$ for
elements of the second, and give the relations between $b,b'$\cite{Ma:exa}. The
abstract theory of braided categories is nevertheless needed to ensure that all
these different braided tensor products are mutually consistent or `coherent'.
Otherwise the idea would be too general and probably intractable.

The appropriate axioms for a $*$-braided group were introduced in \cite{Ma:mec}
for a large class of braided groups and confirmed in \cite{Ma:star} for linear
spaces. We require that our braided group is a $*$-algebra and
\eqn{bg*}{(*\tens *)\circ\und\Delta=\tau\circ\und\Delta\circ*,\quad
\overline{\und\eps(\ )}=\eps\circ*,\quad *\circ \und S=\und S\circ *}
where $\tau$ is the usual transposition. The duality pairing
$\ev=\epsfbox{cup.eps}$ for braided groups $B,C$ is\cite{Ma:introp}
\eqn{bgdual}{\ev(ab,c)=\ev(a,c\Bt)\ev(b,c\Bo),\quad
\ev(b,cd)=\ev(b\Bt,c)\ev(b\Bo,d),\quad \ev(\und Sb,c)=\ev(b,\und Sc)}
for all $a,b\in B$, $c,d\in C$. Here braided coproducts are denoted $\und\Delta
b=b\Bo\tens b\Bt$, etc. (summation implicit). Equivalently (if the braided
antipode is invertible) we can use $\<\ ,\ \>\equiv\ev(\und S^{-1}(\ ),(\ ))$
obeying
\eqn{newbgdual}{\<ab,c\>=\<a,\Psi(b\tens c\Bo),c\Bt\>,\quad
\<b,cd\>=\<b\Bo,\Psi(c\tens b\Bt),d\>,\quad \<\und Sb,c\>=\<b,\und Sc\>}
where we apply $\Psi$ and evaluate its left-hand output with $\<a,\ \>$ and its
right output with $\<\ ,c\Bt\>$, etc. This follows from the
braided-antihomomorphism property\cite{Ma:introp}
\eqn{braant}{\und S\circ\cdot=\cdot\circ\Psi\circ(\und S\tens\und S),\quad
\und\Delta\circ\und S=(\und S\tens\und S)\circ\Psi\circ\und\Delta}
for $\und S$ and similarly (with $\Psi^{-1}$) for $\und S^{-1}$. It is to avoid
the extra braiding that we generally prefer (\ref{bgdual}). In the case of
strict duality where $B=C^\star$ the map $\ev:C^\star\tens C\to\C$
comes also with a coevaluation $\coev=\epsfbox{cap.eps}$, making $C$ a rigid
object in the braided category. When we have $*$-structures then their natural
axioms under duality for quantum groups and braided groups are
\eqn{dual*}{\overline{\<h,a\>}=\<(Sh)^*,a\>,\quad
\overline{\ev(b,c)}=\ev(b^*,c^\star)}
according to \cite{Wor:com} and \cite{Ma:star} respectively. The $*$-operation
$c^\star$ compatible in this way is not necessarily  `unitary'  in a natural
sense  but is often expressible in terms of second $*$-structure $c^*$ which
is. This was already noted for braided linear spaces in \cite{Ma:star} and will
play a role in Sections~4 and~5.

Among the main theorems about braided groups is that if $B$ is any braided
group living in the braided category of $H$-modules ($H$ quasitriangular), then
there is an ordinary quantum group $B\lbiprod H$ called its {\em bosonisation}
and characterised abstractly as such that the ordinary representations of
$B\lbiprod H$ are in 1-1 correspondence with the braided ($H$-covariant)
representations of $B$\cite{Ma:bos}. Explicitly, $B\lbiprod H$ is generated as
an algebra by $H,B$ and has
cross relations, coproduct and antipode
\ceqn{bos}{ hb=(h\o\la b)h\t,\quad \Delta b=b\Bo \CR\ut\tens \CR\uo\la
b\Bt,\quad Sb=(\cu \CR\uo\la\und Sb) S\CR\ut}
for all $h\in H, b\in B$. Here $\cu=(S\CR\ut)\CR\uo$ is an element of $H$. The
algebra is a standard semidirect product by the canonical action $\la$ of $H$
on $B$, while the coalgebra is a semidirect coproduct by the coaction induced
by the universal R-matrix as explained in \cite{Ma:skl}. Likewise if $B$ lives
in the category of $A$-comodules ($A$ dual quasitriangular) then there is an
ordinary quantum group $A\rbiprod B$ called its {\em cobosonisation}, such that
the ordinary corepresentations of $A\rbiprod B$ are in 1-1 correspondence with
the braided ($A$-covariant) corepresentations of $B$\cite{Ma:mec}. Explicitly,
$A\rbiprod B$
is generated as an algebra by $B,A$ and has cross relations, coproduct and
antipode
\eqn{cobos}{ ba=a\o b\bo \CR(b\bt\tens  a\t),\quad \Delta b=b\Bo\bo\tens
b\Bo\bt
b\Bt,\quad S b=(\und S b\bo)Sb\bt}
for all $a\in A, b\in B$. This time the coalgebra is a semidirect coproduct by
the coaction of $A$ on $B$, which we denote $b\bo\tens b\bt$ for the resulting
element of $B\tens A$ (summation implicit).

Finally, we recall a theory of twisting of quasi-quantum groups due in to
Drinfeld\cite{Dri:qua}. A special case of it implies at once that if $H$ is a
quasitriangular Hopf algebra and $\chi\in H\tens H$ a {\em quantum 2-cocycle}
in the sense
\eqn{2-cocycle}{
\chi_{12}(\Delta\tens\id)\chi=\chi_{23}(\id\tens\Delta)\chi,\quad
(\eps\tens\id)\chi=1}
then $H_\chi$ defined by \cite{Dri:qua}
\eqn{hopftwist}{\Delta_\chi=\chi(\Delta\ )\chi^{-1},\quad \eps_\chi=\eps,\quad
\CR_\chi=\chi_{21}\CR\chi^{-1},\quad S_\chi=U(S\ )U^{-1}}
where $U=\chi\uo S\chi\ut$, is again a quasitriangular Hopf algebra. The
cohomological terminology in this context is justified in \cite{Ma:clau}. This
specialisation of Drinfeld's ideas was studied in \cite{GurMa:bra}, and using
the formula for $\Delta U$ given there, it is easy to see that if $H$ is a Hopf
$*$ algebra and $(S\tens S)(\chi^{*\tens *})=\chi_{21}$ then
\eqn{*twist}{*_\chi=(S^{-1}U)((\ )^*)S^{-1}U^{-1}}
makes $H_\chi$ into a Hopf $*$-algebra as well. We will see this in more detail
in the course of a proof in Section~2. The purpose of \cite{GurMa:bra} was to
consider how the corresponding braided groups constructed from $H,H_\chi$ by
transmutation\cite{Ma:bra} are related. We will consider the adjoint of this
question, namely a twisting theory of braided groups such that their
bosonisations are related by quantum group twisting as above. We also consider
how $*$ interacts with the bosonisation construction.

The abstract Sections~2,4 in the paper work over a general field or (with
suitable care) a commutative ring for the purely algebraic results, and over
$\C$ when we discuss $*$. The specific examples in Section~3 based on quantum
enveloping algebras work over formal powerseries $\C[[t]]$ in a deformation
parameter whenever we require directly the quasitriangular structure $\CR$, in
the standard way\cite{Dri}. (The Hopf algebras themselves, their $*$-structures
and their representations do not require this, however, and all work over $\C$
when we use suitable generators, again in the standard way.) All antipodes and
braided antipodes are required for convenience to be invertible, which is in
any case automatic for quasitrianglar and dual-quasitriangular Hopf algebras.

\section{$*$-Structure and Twisting of Bosonisations}

In this section we refine some of the abstract results on the bosonisation
construction\cite{Ma:bos}. This construction can then be used to define
inhomogeneous quantum groups as shown in \cite{Ma:poi}, such as  $q$-Poincar\'e
quantum groups where the momentum sector is a braided covector space (a linear
braided group) living in the braided category of $q$-Lorentz corepresentations
(i.e. it is $q$-Lorentz covariant). Bosonisation consists of adjoining the
$q$-Lorentz sector as a particular semidirect product. Another example of
bosonisation is for a super-Hopf algebra in the category of super-vector spaces
generated by a certain quantum group $\Z_2'$. This is like the Jordan-Wigner
transform which consists in adjoining a degree operator and thereby rendering a
fermionic or super system into a bosonic one. These two settings, Lorentz
covariance and supersymmetry are mathematically unified as cases of one
construction.  As well as the Poincar\'e quantum group, other interesting
applications are to supersymmetry\cite{MacMa:str}\cite{MaPla:uni} and to the
theory of differential calculus on quantum
groups\cite{SchZum:bra}\cite{Dra:bra}. For the moment we proceed quite
generally.

Firstly, it should be perfectly clear that the bosonisation and cobosonisation
constructions (\ref{bos}), (\ref{cobos})  are conceptually dual
to one another. The bosonisation in \cite{Ma:bos} was constructed
diagrammatically by a braided group semidirect product construction
to give a certain braided-Hopf algebra. The Hopf algebra $B\lbiprod H$ contains
$H$ and is arranged so that transmutation from this inclusion\cite{Ma:tra}
reconstructs this braided-Hopf algebra. This step is also diagrammatic. By
turning all diagrams up-side-down one gets the dual construction which gives
(\ref{cobos}). We make a braided semidirect coproduct and use the dual
transmutation theory in \cite{Ma:eul}\cite{Ma:bg}. On the other hand, since not
all readers will be comfortable with the diagrammatic theory, we check the
duality now quite explicitly. Once we have the relevant formulae in detail, we
concentrate with just the bosonisation (\ref{bos}), leaving
the corresponding dual results for the cobosonisation as an easy exercise.

To this end, suppose that $H,A$ are dually paired quantum groups as in
(\ref{hopfdual}) with corresponding quasitriangular and dual-quasitriangular
structure. Let $B$ be a braided group in the category of $H$-modules and $C$ a
braided group in the category of $A$-comodules, which when viewed in the
category of $H$-modules is dually paired to $B$ as in (\ref{bgdual}).

\begin{lemma} With $A,H$ and $B,C$ dually paired as stated, the two ordinary
Hopf algebras $H\lbiprod B$ and $A\rbiprod C$ are dually paired. Between the
various subalgebras the pairing is
\[ \<b,a\>=\eps(b)\eps(a),\quad \<h,a\>={\rm usual},\quad
\<b,c\>=\ev(\und S^{-1}b,c),\quad \<h,c\>=\eps(h)\eps(c)\]
for all $a\in A, h\in H, b\in B, c\in C$.
\end{lemma}
\proof The pairing between subalgebras extends uniquely by the conditions
(\ref{hopfdual}) to the pairing
\eqn{bos-cobos}{\<bh,ac\>=\<\CR\umt h\o,a\>\ev(\und S^{-1}b,\CR\umo h\t\la c)}
between general elements. We check that this is indeed a pairing. \note{Note
the special cases
\[ \<bh,a\>=\<h,a\>\eps(b),\quad \<h,ac\>=\<h,a\>\eps(c),\quad
\<bh,c\>=\ev(\und S^{-1}b,h\la c),\quad \<b,ac\>=\<\CR\umt,a\>\ev(\und
S^{-1}b,\CR\umo\la c)\]
which one checks more easily to be consistent with the pairing conditions.} We
write $\<b,c\>=\ev(\und S^{-1}b,c)$ and need only the properties
(\ref{newbgdual}) mentioned in the preliminaries, so in fact this lemma works
if we suppose $\<\ ,\ \>$ directly without necessarily supposing $\und S^{-1}$.
We denote further distinct copies of $\CR$ by $\CR'$, etc., and consider
general $b,d\in B$, $h,g\in H$, $a\in A$ and $c\in C$. Then
\align{&&\equad \<bh\tens dg, \Delta(ac)\>=\<bh,a\o c\Bo\bo\>\<dg,a\t c\Bo\bt
c\Bt\>\\
&& = \<\CR\umt h\o,a\o\>\<\CR'\umt g\o,a\t c\Bo\bt\>\<b,\CR\umo h\t\la
c\Bo\bo\>\<d,\CR'\umo g\t\la c\Bt\>\\
&&=\<\CR\umt h\th\CR'\umt\o g\o,a\>\<b,\CR\umo h\fo \CR'\umt\t g\t\la
c\Bo\>\<h\o\la d,h\t\CR'\umo g\th\la c\Bt\>\\
&&=\<\CR'\umt h\t g\o,a\>\<b,\CR\umt \CR'\umo\o h\th g\t\la c\Bo\>\<h\o\la d,
\CR\umo \CR'\umo\t h\fo g\th\la c\Bt\>\\
&&=\<\CR'\umt h\t g\o,a\>\<b,\CR\umt\la(\CR'\umo h\th g\t\la c)\Bo\>\<h\o\la d,
 \CR\umo\la(\CR'\umo h\th g\t\la c)\Bt\>\\
&&=\<\CR\umt h\t g\o,a\>\<b(h\o\la d),\CR\umo h\th g\t\la c\>\\
&&=\<b(h\o\la d) h\t g,ac\>=\<bhdg,ac\>}
where the first equality is the coproduct from (\ref{cobos}), the second
evaluates (\ref{bos-cobos}), the third uses the Hopf algebra pairing
(\ref{hopfdual}) and also writes evaluation of the coaction of $c\Bo$ as a
further action on it. As the same time we insert an extra $h\o$ acting on $d$,
and $h\t$ on the other input of $\<\ ,\ \>$, knowing that this is trivial since
$\ev$ and  $\<\ ,\ \>$ between the braided groups are $H$-invariant. For the
fifth equality we use $(h\t\tens h\th)\CR'=\CR'(h\th\tens h\t)$ and then the
cocycle properties of $\CR$ from the axioms of a quasitriangular Hopf algebra.
The sixth equality is the covariance of the braided-coproduct of $C$. Finally,
we use the braided-group duality in the form (\ref{newbgdual}) and recognise
the result as the product from (\ref{bos}). This is half of the lemma. For the
other half, we verify on $b\in B$, $h\in H$, $a,d\in A$ and $c,e\in C$,
\align{&&\equad \<\Delta(bh),ac\tens de\>=\<b\Bo\CR\ut h\o,ac\>\<(\CR\uo\la
b\Bt)h\t,de\>\\
&&=\<\CR'\umt\CR\umt\o h\o,a\>\<\CR''\umt h\th,d\>\<b\Bo,\CR'\umo \CR\umt\t
h\t\la c\>\<b\Bt,\CR\umo \CR''\umo h\fo\la e\>\\
&&=\<\CR\umt h\o,a\>\<h\fo \CR''\umt  ,d\>\<b\Bo,\CR'\umt\CR\umo\o h\t\la
c\>\<b\Bt,\CR'\umo\CR\umo\t h\th \CR''\umo  \la e\>\\
&&=\<\CR\umt h\o,a\>\<h\fo \CR''\umt  ,d\>\<b, (\CR\umo\o h\t\la c) (\CR\umo\t
h\th \CR''\umo  \la e)\>\\
&&=\<\CR\umt\o h\o,a\>\< \CR\umt\t h\t\CR\ut  ,d\>\<b, (\CR\umo\o h\th
\CR\uo\la c) (\CR\umo\t h\fo\la e)\>\\
&&=\<\CR\umt h\o,a d\o\>\,\CR\ut,d\t\>\<b, \CR\umo h\t\la((\CR\uo\la c)e)\>\\
&&=\<bh,ad\o c\bo e\>\CR(c\bt\tens d\t)=\<bh,acde\>}
where the first equality is the coproduct from (\ref{bos}), the second
evaluates the (\ref{bos-cobos}) and also moves $\CR\uo$ acting on $b\Bt$ to
$(S\CR\uo)$ (which becomes $\CR\umo$) acting on the other input of the
braided-group pairing (by its $H$-invariance). The third equality uses the
cocycle property of $\CR^{-1}$ coming from the axioms of a quasitriangular Hopf
algebra, and also moves $h\th\tens h\fo$ to the left of $\CR''$ as $h\fo\tens
h\th$ using these axioms. The fourth equality uses the braided group duality
(\ref{newbgdual}). The fifth equality uses the cocycle and other
quasitriangularity axioms for $\CR$ again. The sixth uses the usual duality
(\ref{hopfdual}) covariance of the braided group product. We then use
(\ref{cobos}) to recognise the result. \endproof

This lemma reworks \cite{Ma:mec} where the duality was given explicitly when
$C=B^\star$ rather than $B=C^\star$ as it is here: Both statements are true.
The astute reader will note also that if these constructions are dual and both
$A$ and $H$ are sub-Hopf algebras the both $A,H$ are also mapped projected
onto. Thus both $A\rbiprod C$ and $B\lbiprod H$ are Hopf algebras with
projection in the sense of Radford\cite{Rad:str}. This observation is due to
the author in \cite{Ma:skl}. \endproof

The lemma is useful when making explicit  dualisations between the comodule and
module setting. For example, the coproduct in (\ref{cobos}) can obviously be
considered as a coaction of $A\rbiprod C$ on $C$ from the right. It is how the
quantum Poincar\'e function algebras in \cite{Ma:poi} coact on the space-time
co-ordinate algebra. Dualising means that in the setting above, $B\lbiprod H$
acts on the spacetime co-ordinate algebra.

\begin{cor} In the setting of the last proposition, $C$ is a left $B\lbiprod
H$-module algebra by
\eqn{bosact}{h\la c=c\bo\<h,c\bt\>,\quad b\la c=\CR\umt\la c\Bo \ev(\und
S^{-1}b, \CR\umo\la c\Bt)}
for $h\in H,b\in B$ and $c\in C$. We call this the {\em fundamental
representation} of the bosonised Hopf algebra.
\end{cor}
\begin{figure}
\epsfbox{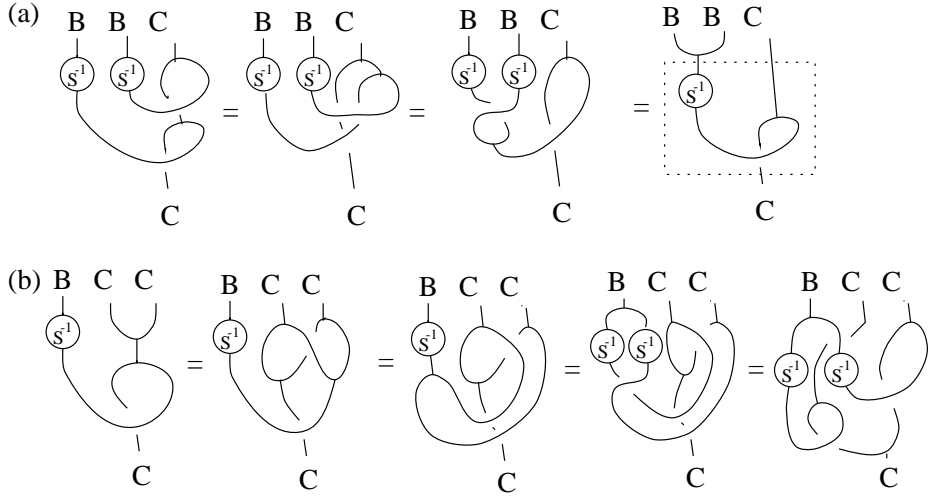}
\caption{Proof that the action of $B$ on $C$ in Corollary~2.2 (in box) is (a)
an action and (b) braided-covariant.}
\end{figure}
\proof Conceptually, the action of $H$ on $C$ is just the action corresponding
to the coaction of $A$ assumed when we said that $B$ was $A$-covariant to begin
with. The action of $B$ in abstract terms is
\eqn{Lreg}{ b\la c=(\ev(\und S^{-1}b,(\
))\tens\id)\circ\Psi^{-1}\circ\und\Delta c}
which has a diagrammatic picture when we write $\Psi$ as a braid crossing and
$\ev=\cup$. This is shown in the box in Figure~1, where we check that it indeed
makes $C$ a braided-module-algebra under its dual $B$. Part (a) checks that it
is an action, using coassociativity of the braided coproduct of $C$, dualising
it to a product of $B$ and the anti-algebra homomorphism property of $\und
S^{-1}$ proven in \cite{Ma:introp}. Part (b) checks that we have a braided
$B$-module algebra in the sense of \cite{Ma:bos}\cite{Ma:introp}, using the
braided bialgebra axiom for $C$, dualising one of its products to a coproduct
of $B$, and the anti-coalgebra homomorphism property of $\und S^{-1}$. The role
of $\ev\circ \und S^{-1}$ can be played directly by a braided group duality
pairing of the type (\ref{newbgdual}) if we prefer. The 1-1 correspondence of
representations in bosonisation theory then means that $C$ becomes an ordinary
module algebra under the bosonisation of $B$ when we adjoin $H$. One gets the
same answer with more work by explicitly evaluating against the coproduct of
$A\rbiprod C$ viewed as a coaction, via the duality pairing from Lemma~2.1.
\endproof

If $B$ is a braided group then its {\em naive opposite coproduct}
$\Psi^{-1}\circ\und\Delta$ makes the algebra of $B$ into a braided group
$B^{\rm cop}$ living not in our original braided category but rather in the
`conjugate' braided category with inverse braiding\cite[Lemma~4.6]{Ma:introp}.
$\und S^{-1}$ becomes its braided
antipode.  In concrete terms it means that the braided group $B^{\rm cop}$ is
no longer properly covariant under $H$ (with the correct induced braiding) but
under this quantum group equipped with $\CR_{21}^{-1}$ instead for its
universal R-matrix. Let us denote the latter by $\bar H$. As a Hopf algebra it
coincides with $H$, but has `conjugate' $\CR$. So we can apply our bosonisation
theorem (\ref{cobos}) and obtain at once a new Hopf algebra $B^{\rm
cop}\lbiprod \bar H$ with coproduct and antipode $\bar\Delta,\bar S$ say. As an
algebra it coincides with $B\lbiprod H$ so $\bar\Delta$ is a second `conjugate'
Hopf algebra structure in this same algebra.

\begin{propos} Let $B\lbiprod H$ be the bosonisation of a braided group $B$.
The second  `conjugate' coproduct and antipode on the same algebra is
\eqn{conjbos}{\bar\Delta b=\CR\umo b\Bt\tens \CR\umt\la b\Bo,\quad \bar
Sb=(\CR\ut\cv^{-1}\la \und S^{-1}b)\CR\uo}
where $\cv=\CR\uo S\CR\ut$, and is twisting equivalent to $(B\lbiprod H)^{\rm
cop}$ by quantum 2-cocycle $\CR^{-1}$. \end{propos}
\proof We use compute from (\ref{bos}) for our new braided group $B^{\rm cop}$
 with opposite coproduct $\Psi^{-1}\circ\Delta b=\CR\umo\la b\Bt\tens
\CR\umt\la b\Bo$ and quantum group $\bar H$ with quasitriangular structure
$\CR_{21}^{-1}$. Then (\ref{bos}) gives bosonised coproduct
\align{\bar\Delta b\equad &&=(\CR\umo\la b\Bt)\CR'\umo\tens
(\CR'\umt\CR\umt)\la b\Bo\\
&&=(\CR\umo\o\la b\Bt)\CR\umo\t\tens \CR\umt\la b\Bo}
using one of the axioms of a quasitriangular structure. We recognise the result
as stated in the proposition when we use the product (\ref{bos}) in $B\lbiprod
H$. The antipode is likewise computed from (\ref{bos}) for our new braided
group $B^{\rm cop}$, which has braided antipode $\und S^{-1}$\cite{Ma:introp}.
We can also compute $\bar\Delta$  further as
\align{\bar\Delta b\equad &&=\CR\umo \CR'\umo\CR\uo b\Bt \tens (\CR\umt\la
b\Bo) \CR'\umt\CR\ut\\
&&=\CR\umo \CR\uo b\Bt \tens (\CR\umt\o\la b\Bo) \CR\umt\t\CR\ut\\
&&=\CR\umo \CR\uo b\Bt \tens \CR\umt b\Bo \CR\ut\\
&&=\CR\umo (\CR\uo\o \la b\Bt)\CR\uo\t \tens \CR\umt b\Bo \CR\ut\\
&&=\CR\umo (\CR\uo\la b\Bt)\CR'\uo \tens \CR\umt b\Bo \CR\ut\CR'\ut}
which we recognise as $\CR^{-1}(\tau\circ\Delta b)\CR$ in view of the coproduct
from (\ref{bos}). The first equality inserts $\CR^{-1}\CR$ into our previous
result for $\bar\Delta b$. The second uses quasitriangularity of $H$, the third
uses the relations in (\ref{bos}) on the right hand factor. We then use the
relations (\ref{bos}) on the left hand factor for the fourth equality and
quasitriangularity again for the fifth. Note that $\bar\Delta h=\Delta h$ since
$H$ is a sub-Hopf algebra, which also equals $\CR^{-1}(\tau\circ\Delta h)\CR$
by quasitriangularity. The quasitriangularity axioms imply that $\CR^{-1}$ is a
2-cocycle for $H^{\rm cop}$ in the sense of (\ref{2-cocycle}), but since $H$ is
a sub-Hopf algebra of $B\lbiprod H$, we
can also view it as a 2-cocycle for $(B\lbiprod H)^{\rm cop}$. We see that the
second `conjugate' Hopf algebra structure on $B\lbiprod H$ is the twisting of
$(B\lbiprod H)^{\rm cop}$ by this cocycle. The antipodes are also twisted one
into the other, since they are determined by the coproducts. \endproof

Both the twisting and the `conjugate' point of view on this second coproduct
are useful.

\begin{corol} In the dual pairing setting of Lemma~2.1,  $C$ is a
module-algebra under $B\lbiprod H$ with its second `conjugate' Hopf algebra
structure, via
\eqn{conjbosact}{h\bar\la c=c\bo\<h,c\bt\>,\quad b\bar\la c=\ev( b,c\Bo)c\Bt}
for $h\in H, b\in B$ and $c\in C$. Moreover, this {\em conjugate fundamental
representation} is isomorphic to the representation in Corollary~2.2 via the
braided antipode of $C$ as intertwiner,
\[ \und S(x\bar\la c)=x\la\und S c ,\qquad\quad\forall x\in B\lbiprod H.\]
\end{corol}
\proof We deduce this without computation via `braided crossing reversal
symmetry', by applying Corollary~2.2 to $B^{\rm cop}$ in the category with
inverse braiding. The role of $C$ is now played by $C^{\rm op}$ with opposite
product $\cdot\circ\Psi^{-1}$. The braided antipode of $C$ is\cite{Ma:introp}
an isomorphism $\und S: C^{\rm op}\to C^{\rm cop}$ of braided groups and we use
this now to refer our action to $C^{\rm cop}$. Then (\ref{Lreg}) uses the
opposite coproduct to the coproduct $\Psi^{-1}\circ\und\Delta$ which (in our
category with reversed braiding) is
$\Psi\circ\Psi^{-1}\circ\und\Delta=\und\Delta$. The $\und S^{-1}$ in
(\ref{Lreg}) becomes the inverse $\und S$ of the antipode of $B$ and is
absorbed in the above isomorphism. In fact, the resulting action is exactly the
left-translation used in defining braided-differentiation in \cite{Ma:fre}, and
we know directly from there that it makes $C$ a braided  $B^{\rm cop}$-module
algebra in the category with reversed braiding (it has the inverse braiding in
its Leibniz rule). It then bosonises to an action of $B\lbiprod H$ by adjoining
the action of $H$. Moreover, we can apply a further $\und S$ to $C^{\rm cop}$
and then our action of $B$ becomes on its image the representation in
Corollary~2.2. Thus
\[ \epsfbox{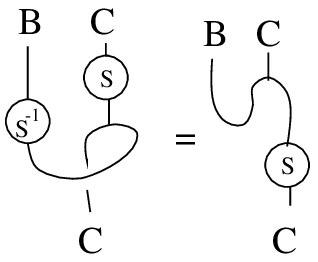}\]
using (\ref{braant}). The action of $H$ in the two cases is the same, namely
the one by which $C$ lives in the category of $H$-modules, and $\und S$ already
intertwines this part of the action because all braided group maps are
morphisms in the category. \endproof

Next we suppose that $H$ is a Hopf $*$-algebra and $B$ is a $*$-braided group,
and that $H$ acts on $B$ `unitarily' in the standard sense
\eqn{act*}{(h\la b)^*=(Sh)^*\la b^*.}
Then the usual theory of Hopf algebra semidirect products ensures that
$B\lbiprod H$ is a $*$-algebra. See \cite{Ma:mec}, where this question was
considered specifically for bosonisations. So $B\lbiprod H$ is certainly a
$*$-algebra.

\begin{lemma} If $H$ is a real-quasitriangular Hopf $*$-algebra and acts
unitarily in the sense (\ref{act*})  on a $*$-braided group $B$ in its category
of representations, then $*$ intertwines the original coproduct $\Delta$ of
$B\lbiprod H$ and $\bar\Delta$, i.e. $(*\tens *)\circ\Delta\circ *=\bar\Delta$.
Likewise, $*\circ S\circ *=\bar S^{-1}$.
\end{lemma}
\proof We compute
\[(*\tens *)\circ\Delta b=\CR\ut^*b\Bo{}^*\tens (\CR\uo\la b\Bt)^*=\CR\umo
b\Bo{}^*\tens\CR\umt\la(b\Bt{}^*)=\bar\Delta(b^*)\]
where the second equality is our reality and unitarity assumption and the third
is the $*$-axiom (\ref{hopf*}) for braided groups. The mapping over of the
antipodes under $*$ as stated is then uniquely determined by the mapping over
of the coproducts. \endproof

Note that the content here is not that $(*\tens *)\circ\Delta\circ *$ is a Hopf
algebra (it is just the Hopf algebra with opposite product, mapped over by $*$)
but that it coincides with $\bar\Delta$ constructed by the `conjugate'
bosonisation. This manifests the deep connection between $*$ and braiding in
the sense of a combined conjugation-braid reversal symmetry. The lemma also
tells us what kind of properties for $*$ to expect for quantum groups
$B\lbiprod H$ obtained by bosonisation. We see that the coproducts $\Delta$ and
$\bar\Delta$ coincide on the quantum group part $H$, where they are both its
usual coproduct. But on the braided group part $B$ they are more
like opposite (transposed) coproducts and indeed become that when $\CR=1$. This
is how our hybrid quantum group interpolates between the axioms (\ref{bg*}) for
a braided group (with a transposition $\tau$) and (\ref{hopf*}) for a usual
quantum group (without $\tau$). Putting this together with the twisting
characterisation of $\bar\Delta$ above, we are motivated to define:

\begin{defin} A quasi-$*$ Hopf algebra is a Hopf algebra which is a $*$-algebra
such that
\cmath{ (*\tens *)\circ\Delta\circ *=\CR^{-1}(\tau\circ\Delta\ )\CR,\quad
\overline{\eps(\ )}=\eps\circ * \\
(\id\tens\Delta)\CR=\CR_{13}\CR_{12},\quad
(\Delta\tens\id)\CR=\CR_{13}\CR_{23},\quad \CR^{*\tens *}=\CR_{21}}
for an invertible element $\CR$ of $H\tens H$.\end{defin}

We will study such objects further in Section~4. One can also consider
something more general where $\CR$ is a cocycle rather than like a
quasitriangular structure. The above definition is stronger but is the one that
applies to our bosonisations. From Proposition~2.3 and Lemma~2.5 we have
clearly:

\begin{corol}  If $B$ is a $*$-braided group acted upon unitarily as in
(\ref{act*}) by a real-quasitriangular $*$-quantum group $H$, then its
bosonisation $B\lbiprod H$ is a quasi-$*$ Hopf algebra with $B,H$ as sub
$*$-algebras.
\end{corol}

Finally, we prove a related and somewhat harder theorem which we will also use
in a later section (as a quantum Wick rotation between Euclidean and Minkowski
Poincar\'e groups.) Namely, we consider how the bosonisation construction
responds to twisting under a cocycle. It is clear that if $B$ is an $H$-module
algebra (an algebra in the category of $H$-modules) and we twist $H$ by a
2-cocycle $\chi$ as in (\ref{hopftwist}) then we must also `twist' the algebra
of $B$ in a certain way if we want it to remain covariant. Likewise, if we have
a coalgebra in the category then we have to `twist' that too if we want to stay
in the category of $H$-modules.

\begin{theorem} If $B$ is a braided group in the category of $H$-modules, and
$\chi$ a 2-cocycle for $H$, then $B_\chi$ defined by
\eqn{bgtwist}{ b\cdot_\chi c=\cdot\circ\chi^{-1}\la(b\tens c),\quad
\und\Delta_\chi b=\chi\la\und\Delta b, \quad \und\eps_\chi=\und\eps,\quad \und
S_\chi=\und S}
is a braided group in the category of $H_\chi$-modules. If $B$ is only a
braided-(bi)algebra then so is $B_\chi$.
\end{theorem}
\proof As explained in \cite{Ma:euc}, for example, we know that if $B$ is a
(co)algebra covariant under $H$ then $B_\chi$ is a (co)algebra covariant under
$H_\chi$. We have to check that this twisted algebra nd coalgebra still fit
together to form a braided group in the braided category of $H_\chi$-modules.
Thus
\align{&&\equad \und\Delta_\chi (b\cdot_\chi c)=\und\Delta_\chi((\chi\umo\la
b)(\chi\umt\la c))\\
&&=\chi\uo\la\left((\chi\umo\la b)\Bo \CR\ut\la(\chi\umt\la
c)\Bo\right)\tens\chi\ut\la\left((\CR\uo\la(\chi\umo\la b)\Bt)(\chi\umt\la
c)\Bt\right)\\
&&=(\chi\uo\o\chi\umo\o\la b\Bo)(\chi\uo\t\CR\ut\chi\umt\o\la c\Bo)\tens
(\chi\ut\o\CR\uo\chi\umo\t\la b\Bt)(\chi\ut\t\chi\umt\t\la c\Bt)\\
&&\equad(\und\Delta_\chi b)\cdot_\chi(\und\Delta_\chi c)=(\chi\umo\chi\uo\la
b\Bo)(\chi\umt\CR\ut_\chi\chi'\uo\la c\Bo)\tens (\chi'\umo\CR\uo_\chi\chi\ut\la
b\Bt)(\chi'\umt\chi'\ut\la c\Bt)}
using the definitions of the product and coproduct of $B_\chi$. In the second
equality we use that $B$ itself is a braided group in the category with
braiding defined via $\CR$, and for the third we use covariance of its product
under $H$. We seek equality for all $b,c\in B$ with the lower expression, which
is the
braided tensor product of $\und\Delta_\chi$ applied to $b,c$ in the category
with braiding defined by $\CR_\chi=\chi_{21}\CR\chi^{-1}$. Equality holds in
view of the identity
\align{&&\equad ((\Delta\tens\Delta)\chi)\CR_{32}\Delta_{H\tens H}\chi^{-1}
=\chi^{-1}_{12}\chi^{-1}_{34}\chi_{23}\CR_{32}\chi^{-1}_{32}\chi_{13}\chi_{24}}
which follows from repeated use of the cocycle condition (\ref{2-cocycle}) and
the quasitriangularity of $\CR$. It is clear that the unit and counit are not
affected by the twisting. Hence we have a braided-bialgebra $B_\chi$. If $B$
has a braided antipode then it is clear that the same map provides a
braided-antipode for $B_\chi$. This is because $\chi$ acts when making the
coproduct  of $B_\chi$ and $\chi^{-1}$ acts when making the product, and $\und
S$ is an intertwiner for the action, so that $\chi^{-1},\chi$ cancel. \endproof

This theorem fits together with the twisting (\ref{hopftwist}) of quantum
groups to tell us that the process of twisting and the process of bosonisation
commute. In categorical terms the reason is that the category of $H_\chi$
covariant representations of $B_\chi$ is equivalent to the category of
$H$-covariant representations of $B$, the equivalence respecting tensor
products up to $\chi$. But the first category is isomorphic to the category of
$B_\chi\lbiprod H_\chi$-modules and the second to that of $B\lbiprod
H$-modules. These are therefore equivalent up to $\chi$. This means by
Tannaka-Krein arguments\cite{Ma:qua} that these two bosonisations are
twisting-equivalent as Hopf algebras.

\begin{propos} In the setting of the preceding theorem, we have
\[ B_\chi\lbiprod H_\chi\isom (B\lbiprod H)_\chi\]
where on the right we view $\chi$ as a 2-cocycle for $B\lbiprod H$. I.e., the
bosonisation of the twisted braided group is the twisted quantum group of its
bosonisation.
\end{propos}
\proof The categorical argument sketched  makes this a corollary of the
preceding theorem. Here we show directly that the required isomorphism is
provided by $\theta:(B\lbiprod H)_\chi\to B_\chi\lbiprod H_\chi$, where
$\theta(bh)=(\chi\uo\la b)\chi\ut h$. This is the identity on the $H_\chi$
sub-Hopf algebra, as it should be since both sides contain this. In addition
\align{&&\equad \theta(bc)=(\chi\uo\la(bc))\chi\ut=(\chi\uo\o\la
b)(\chi\uo\t\la c)\chi\ut=(\chi'\uo\chi\uo\la b)(\chi'\ut\chi''\uo\chi\ut\o\la
c)\chi''\ut\chi\ut\t\\
&&=(\chi\uo\la b)\cdot_\chi(\chi''\uo\chi\ut\o\la
c)\chi''\ut\chi\ut\t=(\chi\uo\la b) \chi\ut(\chi'\uo\la
c)\chi'\ut=\theta(b)\theta(c)}
where we use the $H$-covariance of $B$ for the second equality, the 2-cocycle
condition (\ref{2-cocycle})
for $\chi$ for the third, the definition of the product of $B_\chi$ for the
fourth, and finally the cross relations in $B_\chi\lbiprod H_\chi$ from
(\ref{bos}). We also verify the cross relations as
\[ \theta(h)\theta(b)=h(\chi\uo\la b)\chi\ut=(\chi\uo h\o\la b)\chi\ut
h\t=\theta(h\o\la b)\theta(h\t)\]
where we use the relations in $B_\chi\lbiprod H_\chi$ again. Hence $\theta$ is
an algebra homomorphism.

For the coproduct of $B_\chi\lbiprod H_\chi$ we compute from (\ref{bos})
\align{&&\equad \Delta(\theta(b))=(\chi'\uo\la(\chi\uo\la
b)\Bo)\CR_\chi\ut\chi''\uo\chi\ut\o\chi\umo\tens(\CR_\chi\uo\chi'\ut\la
(\chi\uo\la b)\Bt)\chi''\ut\chi\ut\t\chi\umt\\
&&=(\chi'\uo\chi\uo\o\la
b\Bo)\CR_\chi\ut\chi''\uo\chi\ut\o\chi\umo\tens(\CR_\chi\uo\chi'\ut\chi\uo\t\la
b\Bt)\chi''\ut\chi\ut\t\chi\umt\\
&&\equad \chi((\theta\tens\theta)\circ\Delta b)\chi^{-1}=\chi'\uo(\chi\uo\la
b\Bo)\chi\ut\CR\ut\chi'\umo\tens \chi'\ut(\chi''\uo\CR\uo\la
b\Bt)\chi''\ut\chi'\umt\\
&&=(\chi\uo\chi'\uo\o\la b\Bo)\chi\ut\chi'\uo\t\CR\ut\chi\umo\tens
(\chi''\uo\chi'\ut\o\CR\uo\la b\Bt)\chi''\ut\chi'\ut\t\chi\umt}
where we used the twisted braided coproduct of $B_\chi$ and the twisted
quasitriangular structure of $H_\chi$ for the first equality, and
$H$-covariance of the braided coproduct of $B$ for the second.
We seek equality with the lower expression, computed using the relations of
$B_\chi\lbiprod H_\chi$
obtained from (\ref{bos}). Equality holds for all $b\in B$ in view of the
identity
\align{&&\equad
\chi_{23}\CR_{32}\chi^{-1}_{32}\chi_{13}\chi_{24}(\Delta_{H\tens H}\chi)=
\chi_{12}\chi_{34}((\Delta\tens\Delta)\chi)\CR_{32}}
as in the proof of Theorem~2.8. That the unit and counit map over correctly is
immediate, after which it follows that the antipodes also map over. Hence
$\theta$ is an isomorphism of Hopf algebras. \endproof

We will use this in Section~3.2. For completeness, we also discuss the
interaction of $*$ with twisting.

\begin{propos} If $H$ is a Hopf $*$-algebra and $\chi$ is a 2-cocycle for it
which is {\em real} in the sense $(S\tens S)(\chi^{*\tens *})=\chi_{21}$, then

(i) $H_\chi$ is a Hopf $*$-algebra. If $H$ is (anti)real-quasitriangular then
so is $H_\chi$.

(ii) If $B$ is a $*$-braided-group in the category of $H$-modules with $H$
real-quasitriangular and the action `unitary' in the sense (\ref{act*}), then
$B_\chi$ is a $*$-braided group in the category of $H_\chi$-modules.

(iii) In this case, $B_\chi\lbiprod H_\chi$ is a quasi-* Hopf algebra by
Corollary~2.7.

The required $*$-structures on $H_\chi$ and $B_\chi$ respectively, are
\[ *_\chi(h)=(S^{-1}U)h^*S^{-1}U^{-1},\quad *_\chi(b)=(S^{-1}U)S^{-2}U^{-1}\la
b^*\]
\end{propos}
\proof The first part belongs to the theory of twisting of quantum groups,
mentioned in the Preliminaries. Since it does not seem to be discussed
previously, we include a proof. From the stated `reality' assumption for $\chi$
we see that $U^*=S^{-2}U$ and hence that $S^{-1}U$ is real. This implies that
$(*_\chi)^2=\id$, making $H_\chi$ into a $*$-algebra. Moreover, from the form
of $S_\chi$ in (\ref{hopftwist}) we see that
$S_\chi\circ*_\chi(h)=U(S((S^{-1}U)h^* S^{-1}U^{-1}))U^{-1}=S(h^*)$ so that
$(S_\chi\circ*_\chi)^2=\id$ as required. More non-trivial is the coproduct.
However, it was shown in \cite[Lemma~2.2]{GurMa:bra} that
\eqn{deltaU}{ \Delta U^{-1}=(S\tens S)(\chi_{21})(U^{-1}\tens U^{-1})\chi}
from which we conclude that $\Delta (S^{-1} U^{-1})=\chi^{*\tens
*}(S^{-1}U^{-1}\tens S^{-1}U^{-1})\chi$ under our reality assumption for
$\chi$. This implies at once
\align{&&\equad (*_\chi\tens *_\chi)(\Delta_\chi h)=(S^{-1}U)\chi\umo^* h^*\o
\chi\uo{}^*S^{-1}U^{-1}\tens (S^{-1}U)\chi\umt{}^*
h^*\t\chi\ut{}^*S^{-1}U^{-1}\\
&&\equad =\chi\uo(S^{-1}U)\o h^*\o
(S^{-1}U^{-1})\o\chi\umo\tens\chi\ut(S^{-1}U)\t h^*\ut
(S^{-1}U^{-1})\t\chi\umt=\Delta_\chi\circ*_\chi(h)}
as required. Finally, we check that if $H$ is real-quasitriangular, then
\align{&&\equad (*_\chi\tens *_\chi)(\CR_\chi)=(S^{-1}U\tens
S^{-1}U)\chi^{-1}{}^{*\tens *}\CR_{21}\chi_{21}^{*\tens *}(S^{-1}U^{-1}\tens
S^{-1}U^{-1})\\
&&=\chi(\Delta S^{-1}U) \CR_{21}(\tau\circ\Delta
S^{-1}U^{-1})\chi^{-1}_{21}=(\CR_\chi)_{21}}
using the result on $\Delta S^{-1}U^{-1}$ and the quasitriangularity
assumption. Similarly if $\CR$ is anti-real.

For the second part, we suppose $B$ is a $*$-algebra in the category of modules
of a Hopf $*$-algebra $H$, and that the latter acts as in (\ref{act*}). We
define $*_\chi$ on $B$ as stated. Then
$(*_\chi)^2(b)=(S^{-1}U)S^{-2}U^{-1}\la\left((S^{-1}U)S^{-2}U^{-1}\la
b^*\right)^*=b$ as required, using (\ref{act*}). Moreover, if we write for
brevity $\gamma\equiv(S^{-1}U)S^{-2}U^{-1}$ then our key property
(\ref{deltaU}) tells us that
\eqn{deltagamma}{ \Delta\gamma=\chi^{-1}(\gamma\tens\gamma)(S^{-2}\tens
S^{-2})(\chi).}
Then when we consider the algebra $B_\chi$ with the twisted product
$\cdot_\chi$ as in Theorem~2.8, we will have
\align{&&\equad *_\chi((b\cdot_\chi c))=\gamma\la ((\chi\umt\la
b)^*(\chi\umo\la c)^*)\\
&&=(\gamma\o (S\chi\umt)^*\la c^*)(\gamma\t (S\chi\umo)^*\la
b^*)=(\chi\umo\gamma\la c^*)(\chi\umt\gamma\la
b^*)=(*_\chi(c))\cdot_\chi(*_\chi(b))}
using our reality assumption on $\chi$ and (\ref{deltagamma}). Hence $B_\chi$
is a $*$-algebra under $*_\chi$. Likewise, suppose that $B$ is an
anti-$*$-coalgebra in the sense of the coproduct axiom in (\ref{bg*}) and the
action, and covariant under the Hopf $*$-algebra $H$ as before. Then the
twisted coproduct $\und\Delta_\chi$ as in Theorem~2.8 obeys
\align{&&\equad (*_\chi\tens *_\chi)\und\Delta_\chi(b)=\gamma\la(\chi\uo\la
b\Bo)^*\tens \gamma\la(\chi\ut\la b\Bt)^*\\
&&=\gamma\o (S\chi\umt)^*\la b^*\Bt\tens \gamma\t (S\chi\umo)^*\la
b^*\Bo=\chi\umt\gamma\t\la b^*\Bt\tens\chi\umo\gamma\o\la
b^*\Bo=\tau\circ\und\Delta_\chi\circ *_\chi(b)}
as required.

So in particular, if $B$ is a $*$-braided group then so is $B_\chi$. Moreover,
we know from Theorem~2.8 that $H_\chi$ acts on it. We check that its action
obeys the condition (\ref{act*}). Thus
\[ *_\chi(h\la b)=\gamma(Sh)^*\la
b^*=S^{-1}(U^{-1}(S^{-1}U)h^*(S^{-1}U^{-1})U)\la
*_\chi(b)=S_\chi^{-1}(*_\chi(h))\la *_\chi(b)\]
as required. Since we have also seen that $H_\chi$ is real-quasitriangular when
$H$ is, we are in the situation of Corollary~2.7 and conclude that the
bosonisation $B_\chi\lbiprod H_\chi$ is a quasi-* Hopf algebra for the
$*$-structures as stated. It is related to $B\lbiprod H$ by a theory of
twisting of quasi-* Hopf algebras. \endproof

\section{Poincar\'e quantum enveloping algebras of braided linear spaces}

In this section we specialise our preceding results to the case where $B$ is a
braided vector space $V(R',R)$ as introduced in \cite{Ma:poi}. The data are two
matrices $R',R\in M_n\tens M_n$ obeying the equations
\ceqn{R'R}{R'_{12}R_{13}R_{23}=R_{23}R_{13}R'_{12},\quad
R_{12}R_{13}R'_{23}=R'_{23}R_{13}R_{12}\\
R_{12}R_{13}R_{23}=R_{23}R_{13}R_{12},\quad (PR+1)(PR'-1)=0}
where $P$ is the usual permutation matrix. Then $V$ is defined with generators
$\vecp={p^i}$ and the relations $\vecp_1\vecp_2=R'\vecp_2\vecp_1$, which is
{\em not} the usual Zamolodchikov or exchange algebra since we do not assume
that $R'$ obeys the QYBE (though it often does in practice). Rather, the QYBE
applies to $R$, which we use in defining the braiding
$\Psi(\vecp_1\tens\vecp_2)=R\vecp_2\tens\vecp_1$.
We take braided coproduct $\und\Delta\vecp=\vecp\tens 1+1\tens\vecp$. This can
be expressed as the addition of braided co-ordinates, see
\cite{Ma:poi},\cite{Ma:fre} where the construction was introduced and applied,
respectively.

We assume further that the matrix $R$ is {\em regular} in the
sense\cite{Ma:lin} that we can build from the usual quantum matrix bialgebra
$A(\lambda R)$\cite{FRT:lie} a quantum group $A$ by adding further relations,
such that the canonical (dual) quasitriangular structure on $A(\lambda R)$
introduced in \cite{Ma:qua} descends to $A$. Here $\lambda$ is a constant which
does not enter into the relations of $A(\lambda R)$ but does affect its
dual-quasitriangular structure. It is the {\em quantum group normalisation
constant} introduced in this way in \cite{Ma:lin}. We also make a covariance
assumption\cite{Ma:poi} that the matrix $R'$ is compatible with the quantum
group $A$ in the sense $R'\vect_1\vect_2=\vect_2\vect_1R'$. This is generally
true, for example if $PR'$ is a function of $PR$.

Finally, we assume  that there is a quasitriangular Hopf algebra $H$ dually
paired with $A$. In this case, define the $H$-valued matrices
\eqn{lpm}{ \vecl^+=\<\id\tens \vect,\CR\>,\quad
\vecl^-=\<\vect\tens\id,\CR^{-1}\>.}
which necessarily obey (among other relations) the quadratic relations in
\cite{FRT:lie}. We assume that the elements of these matrices generate $H$  at
least over formal powerseries a deformation parameter. All of the above
conditions are satisfied of course for the standard deformations $U_q(\cg)$
associated to a complex semisimple Lie algebra, but it is necessary for us not
to be tied to this case to cover other interesting examples as well.  This
method to obtain $\vecl^\pm$ from $\CR$ was used in  \cite{Ma:dua} for
$U_q(su_3)$.

This describes all the data for the bosonisation construction, in R-matrix
form. The Poincar\'e enveloping algebras were then constructed under these
assumptions in \cite{Ma:poi} by bosonisation (\ref{bos}), as well as the their
dual quantum groups (which we do not discuss explicitly) by cobosonisation
(\ref{cobos}). To $H$ we have to add a central primitive generator $\xi$ with
quasitriangular structure $\lambda^{-\xi\tens\xi}$ which we multiply into the
quasitriangular structure $\CR$ of $H$ above. This is the quantum group
$\widetilde{H}$ which we actually use. The braided vectors $V(R',R)$ live in
its category of modules by the action
\eqn{lpmcov}{\vecl^+_1\la\vecp_2=\lambda^{-1}R^{-1}_{21}\vecp_2,\quad
\vecl^-_1\la\vecp_2=\lambda R\vecp_2,\quad
\lambda^\xi\la\vecp=\lambda^{-1}\vecp.}
In \cite{Ma:poi} we emphasised the right coaction of the dilatonic extension of
$A$ (the $p^i$ transform as a quantum vector); the above is nothing more than
an evaluation of $\vecl^\pm,\xi$ against that coaction. The {\em inhomogeneous
quantum group} $V(R',R)\lbiprod \widetilde{H}$ of the braided linear space is
then constructed from (\ref{bos}) as generated by $\vecp,\vecl^\pm,\xi$ and
with cross relations, coproduct and antipode\cite{Ma:poi}
\ceqn{vecpoienv}{ \vecl_1^+
\vecp_2=\lambda^{-1}R_{21}^{-1}\vecp_2\vecl^+_1,\quad
\vecl_1^-\vecp_2=\lambda R\vecp_2\vecl^-_1,\quad
\lambda^\xi\vecp=\lambda^{-1}\vecp\lambda^{\xi}\\
\Delta\vecp=\vecp\tens 1+\lambda^\xi\vecl^-\tens\vecp,\quad \eps \vecp=0,\quad
S\vecp=-\lambda^{-\xi} (S\vecl^-)\vecp.}

It should be clear from \cite{Ma:poi} that the reason why there is only
$\vecl^-$  and not $\vecl^+$ in $\Delta \vecp$ is that we used $\CR$ and not
$\CR_{21}^{-1}$ in (\ref{bos}).

\begin{propos} The Hopf algebra $V(R',R)\lbiprod \widetilde{H}$ in the setting
above has a second `conjugate' coproduct and antipode from Proposition~2.3,
namely
\eqn{conjvec}{ \bar\Delta \vecp=\vecp\tens 1+\lambda^{-\xi}\vecl^+\tens
\vecp,\quad \bar S\vecp=-\lambda^\xi(S\vecl^+)\vecp.}
\end{propos}
\proof The computation is exactly the same in form as in \cite{Ma:poi} and
follows at once from (\ref{lpm}), (\ref{lpmcov})  and (\ref{conjbos}). For
example, the `conjugate' coproduct of $\vecp$ is computed as $\vecp\tens
1+\CR\umo\tens \CR\umt\la \vecp$ which yields $\vecl^+$ because the action of
$\CR\umt$ is by evaluation against the coaction $p^i\to p^a St^i{}_a$, which we
compute from (\ref{lpm}) and $(\id\tens S)(\CR^{-1})=\CR$. One may check
directly that the previous coproduct and the new `conjugate' one are related by
conjugation by $\CR^{-1}$, as they must be from Proposition~2.3. \endproof

To proceed further and obtain a quasi-$*$ Hopf algebra, we need to suppose that
$V(R',R)$ is a $*$-braided group. This question was analysed in \cite{Ma:star}
and there are several possibilities. The simplest, which covers $q$-Euclidean
spaces, is to assume that $R$ is of real type I in the sense
$\overline{R}=R^{\dagger\tens\dagger}$ with real quantum group normalisation
constant $\lambda$, and that there is a covariant {\em quantum metric}
$\eta_{ij}$. The construction of such quantum metrics from braided geometry has
been covered in \cite{Ma:eps}. We require a tensor that is quantum group
invariant and obeys various identities with respect to $R,R'$ such as to make
$V(R',R)$ isomorphic the braided covectors $\Vhaj(R',R)$ as braided groups in
the category generated by $R$. The braided covectors are defined in
\cite{Ma:poi} with generators $\{p_i\}$ where the indices are down, and
corresponding relations like the above but $R',R$ acting from left on the
covector generators; we require that $p_i=\eta_{ia}p^a$ effects an isomorphism,
cf\cite{Mey:new}. One may deduce various useful identities, among them (coming
from invariance) the identities\cite{KemMa:alg}
\eqn{metricid}{ \eta_{ia}R^{-1}{}^a{}_j{}^k{}_l=\lambda^2
R^a{}_i{}^k{}_l\eta_{aj},\quad
\eta_{ka}R^i{}_j{}^a{}_l=\lambda^{-2}R^{-1}{}^i{}_j{}^a{}_k\eta_{al}}
which we particularly need below. We use conventions in which $\eta^{ij}$ is
the inverse transpose of $\eta$, and assume further the reality condition
$\overline{\eta_{ij}}=\eta^{ji}$.  Then we can take  $p^i{}^*=p_i$ as explained
in \cite{Ma:star}.

We also need that $H$ is a real-quasitriangular Hopf $*$-algebra. This will be
generally follow from the other conditions as long as $R$ is of real type. For
then the $*$-structure of $A$ can typically be taken in the compact form
$t^i{}_j{}^*=St^j{}_i$ as in \cite{Wor:com}. The dual-quasitriangular structure
of $A$ then necessarily of real type.

\begin{propos} Under the reality assumption on $R$ and given a quantum matric,
the inhomogeneous quantum group $V(R',R)\lbiprod \widetilde{H}$ becomes a
quasi-* Hopf algebra with
\eqn{lpmp*}{ l^\pm{}^i{}_j{}^*=Sl^\mp{}^j{}_i,\quad p^i{}^*=p_i,\quad
\xi^*=\xi}
We regard the $p_i$ as linear combinations of the $p^i$, just with `lowered
indices'.
\end{propos}
\proof The form of $*$ on the $p^i$ is from \cite{Ma:star}. The form on $l^\pm$
is the standard one in the case of $U_q(\cg)$ but is deduced in our more
general setting above from (\ref{lpm}) and the reality type of $\CR$. This was
explained in \cite{Ma:mec}. We also need that the action (\ref{lpmcov}) is
`unitary', which is easily checked from the reality type of $R$. Thus
\align{(l^+{}^i{}_j\la p^k)^*\equad &&
=\overline{\lambda^{-1}R^{-1}{}^k{}_b{}^i{}_j
p^b}=\lambda^{-1}R^{-1}{}^j{}_i{}^b{}_k\eta_{ba}p^a\\
&&=\lambda \eta_{kb} R^j{}_i{}^b{}_ap^a=l^-{}^j{}_i\la \eta_{kb} p^b=(S
l^+{}^i{}_j)^*\la p^i{}^*}
and similarly for the action of $\vecl^-$. We use the equations
(\ref{metricid}).
We then deduce that we have a quasi-$*$ Hopf algebra from Corollary~2.7.
\endproof

Since the $p_i$ are equally good generators in our setting (or if we just want
to work with lower indices in any case) there is nothing stopping us giving the
inhomogeneous quantum group relations in lowered-index form. We can proceed
from the above, using the metric relations (\ref{metricid}) or directly by
bosonisation of $\Vhaj(R',R)$. The results are the same except for a difference
in sign of $\xi$ which can be absorbed by a redefinition. In the first point of
view, the result is
\ceqn{lowvecpoienv}{ \vecp_1\vecp_2=\vecp_2\vecp_1R',\quad \vecl_1^+
\vecp_2=\lambda\vecp_2R_{21}\vecl^+_1,\quad \vecl_1^-\vecp_2=\lambda^{-1}
\vecp_2R^{-1}\vecl^-_1,\quad
\lambda^\xi\vecp=\lambda^{-1}\vecp\lambda^{\xi}\\
\Delta p_i=p_i\tens 1+ \lambda^\xi Sl^-{}^a{}_i\tens p_a,\quad \eps
\vecp=0,\quad Sp_i=-\lambda^{-\xi} (S^2 l^-{}^a{}_i)p_a\\
\bar\Delta p_i=p_i\tens1+\lambda^{-\xi}Sl^+{}^a{}_i\tens
p_a,\quad \bar Sp_i=-\lambda^\xi(S^2 l^+{}^a{}_i)p_a}
where now $\vecp=\{p_i\}$.

The above theory includes, for example, the Euclidean group of motions
$\R_q^n\lbiprod U_q(so_n)$  where $U_q(so_n)$ is the standard $q$-deformed
enveloping algebra, which we take in FRT form\cite{FRT:lie} and Drinfeld's
quasitriangular structure\cite{Dri}. The appropriate $\R_q^n$ are the quantum
planes in \cite{FRT:lie} for suitable $R'$. One of the major results in
\cite{Ma:poi} was that the Poincar\'e quantum groups obtained in this way
automatically (co)act on the spacetime braided covectors $x_i$. As also
explained in \cite{Ma:poi} and developed fully in its sequel \cite{Ma:fre} the
coaction by the braided addition becomes by evaluation an action of the
Poincar\'e enveloping algebra momentum generators $p^i$ on spacetime by braided
differentiation. Since we have two coproducts, we obviously have two such
actions by differentials and `conjugate' differentials. These generate the
fundamental and conjugate fundamental representations from Corollaries~2.2
and~2.4, and are studied further Sections~4 and~5. Let us note that this Hopf
algebra $\R_q^n\lbiprod U_q(so_n)$ and the $q$-euclidean space on which it acts
have recently been studied by rather more explicit means in \cite{Fio:euc}.

\subsection{Spinorial $q$-Euclidean Poincar\'e enveloping algebra}

There are many classes of linear braided groups, much beyond the usual quantum
planes associated to representations of standard quantum group deformations
such as $\R_q^n\lbiprod U_q(so_n)$. In this subsection and the next, we
specialise to a class in which the generators $p_i$ above are replaced by a
matrix of generators $\vecp=\{p^{i_0}{}_{i_1}\}$ say. We are still considering
them as an additive braided group but adopt a notation in which $\vecp$ refers
to a matrix, with the usual notational rules. For example, the matrix
bialgebras $A(R)$ have such a coaddition\cite{Ma:add}, as do rectangular
quantum matrices $A(R:S)$ in \cite{MaMar:glu}, whenever $R,S$ are $q$-Hecke
solutions of the QYBE (so that the corresponding braiding has eigenvalues
$q,-q^{-1}$).

In particular, we focus here on  the `rectangular quantum matrix' algebra $\bar
A(R)=A(R_{21}:R)$ proposed in the $2\times 2$ case as a matrix or `twistor'
description of 4-dimensional $q$-Euclidean space in \cite{Ma:euc}. It has
relations and braid statistics
\eqn{eucmat}{R_{21}\vecp_1\vecp_2=\vecp_2\vecp_1R,\quad
\vecp'_1\vecp_2=R\vecp_2\vecp_1'R.}
under which it forms a braided group with $\und\Delta\vecp=\vecp\tens
1+1\tens\vecp$ as before. Clearly we can write $p^{i_0}{}_{i_1}=p_I$ as a
covector with multiindex $I=(i_0,i_1)$, and in this way write this braided
linear space as $\Vhaj({\bf R}',{\bf R})$ for suitable `multiindex' ${\bf
R}',{\bf R}$ given in \cite{Ma:euc}. So the difference is purely notational.

One can then follow the preceding section with the matrices ${\bf R}',{\bf R}$
and quantum group generators $l^\pm{}^I{}_J$ using the setting there.
Alternatively, which we do in the present section, we can take a slightly
different quantum group as the `background symmetry' with respect to which we
bosonise. Namely, we take in the role of $H$ in the preceding section a quantum
group $H\tens H$, where this time $H$ is a quasitriangular Hopf algebra dual to
a dual-quasitriangular Hopf algebra $A$ obtained from $A(\lambda^\h R)$. We
take the tensor product quasitriangular Hopf algebra structure. This quantum
group is related to the one in the preceding section by the realisation
\eqn{euclreal}{l^\pm{}^I{}_J
=(S^{-1}l^\pm{}^{j_0}{}_{i_0})m^\pm{}^{i_1}{}_{j_1}}
where $\vecl^\pm$ and $\vecm^\pm$ denote the matrix generators of the two
copies of $H$. It is still the case that $\bar A(R)$ lives in the braided
category of a dilatonic extension of $H\tens H$, so we construct this and
proceed directly by bosonisation. This formulation has been explained in
\cite{Ma:euc} and we already know from there that the spacetime-coordinates
become a module algebra under the spacetime rotation. We use the same action
on the lowered-index momentum generators $p_I=p^{i_0}{}_{i_1}$, regarded
now as a matrix, namely \cite{Ma:euc}
\eqn{lmcov}{\vecl^+_1\la \vecp_2=\lambda^{-\h}R^{-1}_{21}\vecp_2,\quad
\vecl^-_1\la\vecp_2=\lambda^\h R\vecp_2,\quad \vecm^+_1\la \vecp_2=\vecp_2
\lambda^\h R_{21},\quad \vecm^-_1\la\vecp_2=\vecp_2 \lambda^{-\h}R^{-1}.}
We add a dilaton $\xi$ as before, with $\lambda$ the quantum group
normalisation constant for ${\bf R}$, which is the square of that for $R$.

\begin{propos} The inhomogeneous quantum group $\bar A(R)\lbiprod
(\widetilde{H\tens H})$ constructed from (\ref{bos}) has cross relations,
coproduct and antipode
\ceqn{eucspinpoienv}{ \vecl^+_1
\vecp_2=\lambda^{-\h}R^{-1}_{21}\vecp_2\vecl^+_1,\quad
\vecl^-_1\vecp_2=\lambda^\h R\vecp_2\vecl^-_1\\
 \vecm^+_1 \vecp_2=\vecp_2\lambda^\h R_{21}\vecm^+_1,\quad
\vecm^-_1\vecp_2=\vecp_2 \lambda^{-\h}R^{-1}\vecm^-_1,\quad
\lambda^\xi\vecp=\lambda^{-1}\vecp\lambda^\xi\\
\Delta \vecp=\vecp\tens 1+\lambda^\xi (\vecl^-(\ )S\vecm^-)\tens\vecp,\quad
\eps \vecp=0,\quad S\vecp=-\lambda^{-\xi}S(\vecl^-(\ )S\vecm^-)\vecp}
where $\vecl^-(\ )S\vecm^-$ has a space for the matrix indices of $\vecp$ to be
inserted. The second `conjugate' coproduct and antipode from Proposition~2.3
are
\eqn{conjeucspin}{\bar\Delta \vecp=\vecp\tens 1+\lambda^{-\xi}\vecl^+(\
)S\vecm^+\tens\vecp,\quad \bar S\vecp=-\lambda^\xi S(\vecl^+(\
)S\vecm^+)\vecp.}
\end{propos}
\proof The semidirect product  (\ref{bos}) gives the cross
relations as before. We read them off from (\ref{lmcov}) because of the matrix
form of the coproducts of $\vecl^\pm,\vecm^\pm$. For the coproduct and
conjugate coproduct we evaluate against the composite $\CR$ (which is the
tensor product of one for each copy). The computation follows the same line as
in \cite{Ma:poi} and Proposition~3.1, for each copy separately, giving the
result. \endproof

This has been announced in \cite{Ma:clau94}. It is dual by Lemma~2.1 to the
spinorial Poincar\'e enveloping algebra computed in \cite{Ma:euc}. Next we
consider $*$-structures. One can take various $*$-structures on $\bar A(R)$
(including one as a Hermitian matrix) but we concentrate here on  `unitary
type' $*$-structures defined according to the twisting theory in \cite{Ma:euc}
by having the same form as on the generators $\vect$ of the quantum group $A$
obtained from $A(\lambda^\h R)$. We suppose this has the explicit form
$t^i{}_j{}^*=\eps_{ai}t^a{}_b\eps^{bj}$
say, where $\eps^{ij}$ is quantum group invariant and $\eps_{ij}$ the
transposed inverse. For $H\tens H$, the $*$-structure we need according to the
theory in \cite{Ma:euc} is not quite the tensor product one, but has an extra
automorphism by $S^{-2}$ in the first copy. We assume as before that $H$ is
real quasitriangular Hopf dual to $A$. For both of these assumptions we assume
that $R$ is of real type I in the sense $\bar R=R^{\dagger\tens\dagger}$ and
that $\lambda^\h$ is real.

\begin{propos} Under the reality assumption on $R$ and given a suitable $\eps$
as above, the inhomogeneous quantum group $\bar A(R)\lbiprod\widetilde{H\tens
H}$ becomes a quasi-$*$ Hopf algebra with
\eqn{lm*}{ l^\pm{}^i{}_j{}^*=S^{-1}l^\mp{}^j{}_i,\quad
m^\pm{}^i{}_j{}^*=Sm^\mp{}^j{}_i,\quad
p^i{}_j{}^*=\eps_{ai}p^a{}_b\eps^{bj},\quad \xi^*=\xi}
\end{propos}
\proof The quantum group $H\tens H$ is real-quasitriangular since each factor
is. The $*$ structure on $\vecm^\pm$ is as in the preceding section. On the
$\vecl^\pm$ we have the extra automorphism $S^{-2}$. We also need the action
(\ref{lmcov}) to be unitary in the sense (\ref{act*}). Thus
\align{(m^+{}^i{}_j\la p^k{}_l)^*&&\equad =(p^k{}_a\lambda^\h
R^a{}_l{}^i{}_j)^*=\lambda^\h R^j{}_i{}^l{}_a  \eps_{bk}p^b{}_c\eps^{ca}=
\lambda^{-\h}p^b{}_c R^{-1}{}^j{}_i{}^c{}_a\eps_{bk}\eps^{al}\\
&&=m^-{}^j{}_i\la\eps_{bk}p^b{}_c\eps^{cl}=S^{-1}(m^+{}^i{}_j{}^*)\la
p^k{}_l{}^*\\
(l^+{}^i{}_j\la p^k{}_l)^*&&\equad =(p^k{}_a\lambda^{-\h}
R^{-1}{}^k{}_a{}^i{}_j)^*=\lambda^{-\h}
R^{-1}{}^j{}_i{}^a{}_k\eps_{ba}p^b{}_c\eps^{cl}
=\lambda^{\h}\widetilde{R^{-1}}{}^j{}_i{}^a{}_b \eps_{ak}p^b{}_c\eps^{cl}\\
&&=S^{-2}l^-{}^j{}_i\la\eps_{bk}p^b{}_c\eps^{cl}=S^{-1}(l^+{}^i{}_j{}^*)\la
p^k{}_l{}^*}
using for $\vecm^+$ the equations (\ref{metricid}) for our invariant tensor
$\eps$, and a variant of them (proven in the same way) for $\vecl^+$. Here
$\widetilde{R^{-1}}$ is the `second inverse' of $R^{-1}$ and governs the action
of $S^{-2}\vecl^-$ deduced from (\ref{lmcov}). Similarly for the
$\vecm^-,\vecl^-$ cases. We then deduce that we have a quasi-* Hopf algebra
from Corollary~2.7. \endproof

We can put general Hecke R-matrices into the above constructions. For the
standard $su_2$ R-matrix we have for the braided linear group the rectangular
quantum matrices $\bar M_q(2)$ which are isomorphic, in this particular case,
to usual $M_q(2)$. In this way we have compatibility with a previous proposal
for a suitable algebra for 4-dimensional Euclidean space in \cite{CWSSW:ten}.
Here $\bar M_q(2)\lbiprod \widetilde{U_q(su_2)\tens U_q(su_2)}$ is a
`spinorial' version of the $q$-deformed Euclidean group of motions.

\subsection{Spinorial $q$-Minkowski Poincar\'e enveloping algebra}

In this section we consider another braided linear space in matrix form
$\vecp=\{p^{i_0}{}_{i_1}\}$, namely the braided matrices $B(R)$ introduced by
the author as a multiplicative braided group in \cite{Ma:exa}. The additive
braided group structure is due to U. Meyer in \cite{Mey:new} and requires that
$R$ is $q$-Hecke, which we assume. The relations and additive braid statistics
are
\eqn{bramat}{ R_{21}\vecp_1R\vecp_2=\vecp_2R_{21}\vecp_1R,\quad R^{-1}\vecp_1'
R\vecp_2=\vecp_2 R_{21}\vecp_1'R}
and we take braided coproduct $\und\Delta\vecp=\vecp\tens 1+1\tens\vecp$. As
before, we can also write $p^{i_0}{}_{i_1}=p_I$ as a braided covector space
$\Vhaj({\bf R}',{\bf R})$ for suitable ${\bf R}',{\bf R}$ given in
\cite{Ma:exa}, \cite{Mey:new} respectively. The equivalence between the
notations is standard after the paper  \cite{Ma:skl}. The algebra relations in
(\ref{bramat}) are of interest in  other contexts
too\cite{FRT:lie}\cite{ResSem:mat}, as explained in \cite{Ma:skl}.

One can then follow the first subsection above with the matrices ${\bf R}',{\bf
R}$ and quantum group generators $l^\pm{}^I{}_J$ using the setting there. This
approach to the $q$-Lorentz group is covered in Meyer's paper \cite{Mey:new}.
Alternatively, which we do, we can follow the `spinorial' point of view and
take for our background quantum group symmetry the quantum group $H\codcross H$
 obtained by twisting the quasitriangular Hopf algebra $H\tens H$ in the
preceding section by the quantum 2-cocycle $\chi=\CR_{23}^{-1}$ as an element
of $(H\tens H)^{\tens 2}$. I.e., $\chi$ is $\CR^{-1}$ but with its first
component living in the copy of $H$ with generators $\vecm^\pm$ and its second
component living in the copy of $H$ with generators $\vecl^\pm$. The coproduct,
antipode and quasitriangular structure are read off from (\ref{hopftwist}).
Note that the use of  this twisted Hopf algebra to describe the Lorentz quantum
group is due to the author in \cite[Sec. IV]{Ma:poi}, where we pointed out for
the first time the isomorphism of previous proposals for the Lorentz quantum
group function algebra  in \cite{CWSSW:lor}\cite{PodWor:def} with the dual of
the `twisted square' in \cite{ResSem:mat}. It has subsequently been reiterated
by other authors. The realisation of this quantum group in terms of the
vectorial picture in the first subsection is cf\cite{Mey:new}
\eqn{Llm}{
l^+{}^I{}_J=(\vecl^-\vecm^+)^{i_1}{}_{j_1}((S_0^{-1}\vecm^+)
(S_0^{-1}\vecl^+))^{j_0}{}_{i_0},\quad
l^-{}^I{}_J=(\vecl^-\vecm^-)^{i_1}{}_{j_1}((S_0^{-1}\vecm^+)
(S_0^{-1}\vecl^-))^{j_0}{}_{i_0}}
where $S_0$ is the usual `matrix inverse' antipode of $H$. The two copies of
$H$ no longer appear as sub-Hopf algebras due to the twisting. It is still the
case that $B(R)$ lives in the category of representations of a dilatonic
extension of $H\codcross H$, and we bosonise with respect to this. The required
action on the spacetime co-ordinates $B(R)$, which we use now on the lowered
momentum generators $p_I=p^{i_0}{}_{i_1}$ has already been given in
\cite{Ma:euc}: by the twisting theory developed there, we use exactly the same
formula (\ref{lmcov})
on the generators, but extending now to our new algebras. We add the dilaton
$\xi$ as before, with $\lambda$ the quantum group normalisation constant of
${\bf R}$, which is again the square of that of $R$.

\begin{propos} The inhomogeneous quantum group $B(R)\lbiprod
(\widetilde{H\codcross H})$ constructed from (\ref{bos}) has cross relations,
coproduct and antipode
\ceqn{spinpoienv}{R_{21}\vecp_1R\vecp_2=\vecp_2R_{21}\vecp_1R,\quad
\vecl^+_1\vecp_2\vecl^-_2=\lambda^{-\h}R^{-1}_{21}\vecp_2\vecl^-_2\vecl^+_1,
\quad \vecl^-_1\vecp_2\vecl^-_2=\lambda^\h R\vecp_2\vecl^-_2\vecl^-_1\\
 \vecm^+_1 \vecp_2\vecl^-_2=\vecp_2\vecl^-_2\lambda^\h R_{21}\vecm^+_1,
\quad\vecm^-_1\vecp_2\vecl^-_2=\vecp_2\vecl^-_2
\lambda^{-\h}R^{-1}\vecm^-_1,\quad
\lambda^\xi\vecp=\lambda^{-1}\vecp\lambda^\xi\\
\Delta \vecp=\vecp\tens 1+\lambda^\xi \vecl^-\vecm^+(\
)(S_0\vecm^-)(S_0\vecl^-)\tens\vecp,\quad \eps \vecp=0\\
 S\vecp=-\lambda^{-\xi}S_0\left(\vecm^+\vecl^-(\
)(S_0\vecl^-)(S_0\vecm^-)\right)\vecp}
where $(\ )$ is a space for the matrix entries of $\vecp$ to be inserted, and
$S_0$ is the usual matrix antipode in either copy of $H$. The second
`conjugate'
coproduct and antipode from Proposition~2.3 are
\eqn{conjspin}{\bar\Delta \vecp=\vecp\tens 1+\lambda^{-\xi}\vecl^+\vecm^+(\
)(S_0\vecm^+)(S_0\vecl^-)\tens\vecp,\quad \bar S\vecp=-\lambda^\xi
S_0(\vecm^+\vecl^+(\ )(S_0\vecl^-)(S_0\vecm^+))\vecp.}
\end{propos}
\proof We compute again from (\ref{bos}). This time we evaluate the coaction
$p^i{}_j\to p^a{}_b (Ss^i{}_a)t^b{}_j$ which underlies (\ref{lmcov}, where
$\vecs,\vect$ are dual to $\vecl^\pm,\vecm^\pm$, against the twisted
quasitriangular structure $\chi_{21}\CR_{H\tens H}\chi^{-1}$, using
(\ref{lpm}). The computation follows the same line as in  Proposition~3.3
except that we use the matrix coproduct of $\vecs,\vect$ in the duality pairing
(\ref{hopfdual}) to evaluate products of $\CR$. \endproof

This was announced in \cite{Ma:clau94}. It is dual via Lemma~2.1 to the
spinorial Poincar\'e function algebra computed in \cite{MaMey:bra}. We are now
in a position to say rather more about its structure.

\begin{propos} The quantum group $B(R)\lbiprod(\widetilde{H\codcross H})$ is
the twisting of the quantum group $\bar A(R)\lbiprod(\widetilde{H\tens H})$
from Proposition~3.3, by the quantum 2-cocycle $\chi\in (H\tens H)^{\tens 2}$
viewed in the latter quantum group.
\end{propos}
\proof We just apply Proposition~2.9. That the algebras are indeed isomorphic
is quite easy to see explicitly: we identify $\vecp$ in Proposition~3.3 with
$\vecp\vecl^-$ in Proposition~3.5. That the coalgebras are then related by
twisting requires rather more work to verify directly. \endproof

Thus the two systems based on $\bar A(R)$ and $B(R)$ are algebraically `gauge
equivalent'\cite{Dri:qua} in the sense of twisting of quantum and braided
groups, so which one chose to work with is primarily a matter of convenience
like a `choice of co-ordinates'. The proposition extends this `quantum Wick
rotation' from \cite{Ma:euc} to the level of the associated `Poincar\'e quantum
groups' in the interpretation there. Next we consider $*$-structures. We
suppose that $R$ is of real type I and $\lambda$ real. A natural $*$ structure
was introduced in \cite{Ma:mec}, namely the Hermitian one. For $H\codcross H$
we take the dual of the $*$-structure on the quantum group $A\dcross A$
introduced in \cite{Ma:poi} in our abstract approach to the $q$-Lorentz group
function algebra. On matrix generators it is $s^i{}_j{}^*=St^j{}_i$ as studied
in \cite{CWSSW:lor}.

\begin{propos} Under the reality assumption on $R$, the inhomogeneous quantum
group $B(R)\lbiprod(\widetilde{H\codcross H})$ becomes a quasi-$*$ Hopf algebra
with
\eqn{mink*}{ l^\pm{}^i{}_j{}^*=Sm^\mp{}^j{}_i,\quad
m^\pm{}^i{}_j{}^*=Sl^\mp{}^j{}_i,\quad p^i{}_j{}^*=p^j{}_i,\quad \xi^*=\xi}
\end{propos}
\proof The  antipode  on $H\codcross H$ is the twisted form $S(h\tens
g)=U(S_0h\tens S_0g)U^{-1}$ from (\ref{hopftwist}) where $U=\CR_{21}$, and the
$*$-structure likewise has the twisted form  $(h\tens g)^*=U(g^*\tens
h^*)U^{-1}$ (this is obvious from \cite{Ma:poi} as the dual of the
$*$-structure on $A\dcross A$ there). We take on $H\tens H$ the `flipped'
$*$-structure $(h\tens g)^*=g^*\tens h^*$, with respect to which our 2-cocycle
$\chi$ is of real-type in the sense needed for (\ref{*twist}). We also have
$S_{H\tens H}^{-1}U=U$ for our particular 2-cocycle. So
$l^\pm{}^i{}_j{}^*=S_0m{}^\mp{}^j{}_i$ in $H\tens H$  twists to $H\codcross H$
by the same conjugation factor as for the antipode $S$, giving the form stated.
We also know from Proposition~2.10 that the twisted $\CR$ will be
real-quasitriangular since it clearly is so on $H\tens H$ before twisting. We
check finally that the action (\ref{lmcov}) is unitary in the sense
(\ref{act*}). Thus
\align{(m^+{}^i{}_j\la p^k{}_l)^*\equad &&=(p^k{}_a\lambda^\h
R^a{}_l{}^i{}_j)^*=\lambda^\h R^j{}_i{}^l{}_a p^a{}_k=S^{-1}(Sl^-{}^j{}_i)\la
p^l{}_k=S^{-1}(m^+{}^i{}_j{}^*)\la p^k{}_l{}^*}
as required. Similarly for the action of $\vecm^-$ and $\vecl^\pm$. We then
conclude that we have a quasi-* Hopf algebra from Corollary~2.7. \endproof

This proposition confirms that the present system based on $B(R)$ differs,
however, by more than just a `change of co-ordinates' from the system based on
$\bar A(R)$ from the previous subsection, because it has quite a different
$*$-structure: even if we refer both systems to the same algebra  by untwisting
the $*$-structure in Proposition~3.7, we will not obtain the previous quasi-*
Hopf algebra in Proposition~3.4. Indeed, it is clear from \cite{Ma:poi} and
from the above proof that  $B(R)\lbiprod (\widetilde{H\codcross H})$ is the
twisting via Proposition~2.10 of $\bar A(R)\lbiprod (\widetilde{H\tens H})$
with the `flip' $*$-structure on $H\tens H$ and (since $(S_{H\tens
H}^{-1}U)S_{H\tens H}^{-2}U^{-1}=1$) the same Hermitian $\vecp$, in contrast to
Section~3.1 where we had essentially a tensor product $*$ structure on $H\tens
H$ and a `unitary' type $*$-structure on $\vecp$. Let us note also that while
it may be useful to untwist in order to make such comparisons, there are good
reasons too to work with the $B(R)$ `co-ordinates' most of the time, such as
its multiplicative braided group structure\cite{Ma:exa} and the remarkable
identification of that with the braided-universal enveloping algebra of a
braided Lie algebra $\CL(R)$ \cite{Ma:lie}.

For the standard $su_2$ R-matrix we obtain the braided matrices $BM_q(2)$ in
\cite{Ma:exa}, isomorphic to an algebra proposed  as $q$-Minkowski space in
\cite{CWSSW:ten} from consideration of the tensor product of two quantum plane.
Then $BM_q(2)\lbiprod(\widetilde{H\codcross H})$ is a `spinorial' version of
the $q$-deformed Poincar\'e enveloping algebra in 4-dimensional $q$-Minkowski
space. Such an algebra has been studied in \cite{OSWZ:def} via explicit
generators and relations, i.e the R-matrix form above and the results about it
are new. The braided matrices $BM_q(2)$ are also isomorphic to a degenerate
form of the Sklyanin algebra\cite{Ma:skl} and to the braided enveloping algebra
of the braided-Lie algebra $\und{gl}_q$\cite{Ma:mex}.

\section{Unitary Representations of Quasi-$*$ Hopf Algebras}

In this section we provide some basic lemmas about quasi-* Hopf algebras  and
their representations. We examine in detail the notion of tensor products of
unitary representations. This leads to a general construction for sesquilinear
forms or `inner products' such that the fundamental and conjugate fundamental
representations of our inhomogeneous quantum groups are mutually adjoint and
hence unitary in our sense. This underpins our remarks about the differential
representation of $q$-Poincar\'e quantum groups in the next section.

\begin{lemma} If $(H,\CR,*)$ is a quasi-* Hopf algebra in the sense of
Definition~2.6 then

(i) $(\eps\tens\id)\CR=1=(\id\tens\eps)\CR$, $(S\tens \id)(\CR)=\CR^{-1}$ and
$(S\tens S)(\CR)=\CR$.

(ii) $\CR$ obeys the QYBE in $H\tens H\tens H$

(iii) $\CR$ is a 2-cocycle for $H$ (or equivalently $\CR^{-1}$ is a 2-cocycle
for $H^{\rm cop}$).

(iv) $*\circ S\circ *=\cu^{-1}(S\ )\cu$ where $\cu^{-1}=\CR\ut S^2\CR\uo$ is
invertible and  $\Delta\cu^{-1}=(\cu^{-1}\tens\cu^{-1})\CR_{21}\CR$.
\end{lemma}
\proof These facts are analogous to similar facts for quasitriangular
Hopf algebras\cite{Dri} but require a little more work, except for (i), for
which the proof is unchanged. For (ii) we compute
\align{&&\equad \CR_{13}\CR_{12}=(\id\tens\Delta)\CR=(*\tens \Delta\circ
*)\CR_{21}=(\CR^{-1}_{23}((\id\tens
\tau\circ\Delta)\CR_{21})\CR_{23})^{*\tens*\tens*}\\
&&=(\CR^{-1}_{23}\CR_{31}\CR_{21}\CR_{23})^{*\tens*\tens*}
=\CR_{32}\CR_{12}\CR_{13}\CR^{-1}_{32}}
which is the QYBE in $H\tens H\tens H$ after suitable renumbering. From (ii) we
deduce (iii) at once in view of the existing assumptions for $\Delta$ on $\CR$.
Part (iv) then follows from part (iii) and the theory of twisting of Hopf
algebras cf\cite{Dri:qua}, at least if $S$ is invertible. We view $\CR^{-1}$ as
a cocycle for $H^{\rm cop}$ and deduce that the twisted coproduct
$\bar\Delta=\CR^{-1}(\tau\circ\Delta(\ ))\CR$ has an antipode $\bar S=U(S^{-1}\
)U^{-1}$ where $U=\CR\umo S^{-1}\CR\umt=(S^2\CR\uo)\CR\ut$ using part (i), and
where (from the twisting theory) this is invertible. We denote $U^{-1}=\cv$,
say. But clearly $*\circ S^{-1}\circ *$ is also the antipode for $\bar
\Delta=(*\tens *)\circ\Delta\circ *$, hence by uniqueness of the antipode we
deduce $*\circ S^{-1}\circ *=\cv^{-1}(S^{-1}\ )\cv$. This inverts to the form
stated, where $\cu=S\cv$. Finally, we combine (\ref{deltaU}) from
\cite{GurMa:bra} with (i) and the reality condition $\CR^{*\tens *}=\CR_{21}$
to deduce
$\tau\circ\Delta \cv=\CR^{-1}_{21}(\cv\tens
\cv)\CR^{-1}=(\cv\tens\cv)\CR^{-1}_{21}\CR^{-1}$. The same form for $\cu$
follows
by applying $S$. \endproof

These are some of the most basic features. We denote $H$ with its second
`conjugate' Hopf algebra structure by $\bar H$ (it is also a quasi-$*$ Hopf
algebra, with a different 2-cocycle). It is clear that a quasi-$*$ Hopf algebra
is quasitriangular {\em iff} it is a usual Hopf $*$-algebra (in which case it
is real-quasitriangular), which is {\em iff} $H=\bar H$. Unlike usual Hopf
$*$-algebras, however, we do not generally have $(S\circ*)^2=\id$. This map
$(S\circ *)^2=S\circ\bar S^{-1}$ remains, however, an interesting algebra
automorphism and at least in some contexts it is natural to ask that it be
inner (e.g. if one wants to build a Tannakian category of representations along
the lines for algebraic groups in \cite{DelMil:tan}.)  Part (iv) of Lemma~4.1
in the form
\eqn{S*}{ (S\circ*)^2=\cv(S^2(\ ))\cv^{-1}}
tells us that this is inner {\em iff} the automorphism $S^2$ is inner. For our
bosonisation examples we have:

\begin{propos} (i) Let $H$ be a quasitriangular Hopf algebra and $\alpha\in H$
such that $\Delta\alpha=(\alpha\tens\alpha) (\CR_{21}\CR)^{-1}$. Then $\alpha$
induces an automorphism   $\theta_{\alpha}:B\to B$, $\theta_\alpha(b)=\alpha\la
\und S^2(b)$ of any braided group $B$ in the category of $H$-modules. Here
$\und S$ denotes its braided-antipode.

(ii)  The square of the antipode of the bosonisation Hopf algebras $B\lbiprod
H$ is
\[ S^2(b)=\theta_\cu(b),\quad S^2(h)=\cu h\cu^{-1}\]
for all  $b\in B, h\in H$. Here $\cu=(S\CR\ut)\CR\uo$  as in \cite{Dri:alm}.
\end{propos}
\proof That $\theta_\alpha$ is always a braided group automorphism is clear
from (\ref{braant}): the $\Psi^2$ from the action of $\und S^2$ is cancelled by
the $(\CR_{21}\CR)^{-1}$ in the coproduct of $\alpha$, which determines its
action on tensor products. This part is an elementary fact about braided
groups. In our bosonisation Hopf algebras we compute from (\ref{bos}) that
\align{&&\equad S^2(b)= S((\cu\CR\uo\la\und S
b)S\CR\ut)=(S^2\CR\ut)S(\cu\CR\uo\la \und S b)\\
&&=\CR\ut(\cu\CR'\uo \CR\uo\cu\la\und S^2 b)S\CR'\ut=\CR\ut\cu\CR\uo\cu\la \und
S^2 b=\theta_\cu(b)}
where we use the definition of the antipode in $B\lbiprod H$, the standard fact
that $\cu(\ )\cu^{-1}=S^2$ in $H$ and then the relations in (\ref{bos}).
\endproof

Let us note also that many properties of quasi-$*$ Hopf algebras depend only on
the feature that $(*\tens *)\circ\Delta\circ *$ is twisting equivalent to
$\tau\circ\Delta$. I.e. we can demand only the 2-cocycle property (iii) in
Lemma~4.1 in place of the more restrictive axioms for $(\Delta\tens\id)\CR$ and
$(\id\tens\Delta)\CR$ in Definition~2.6. It is natural to call this a {\em
cocycle}-$*$ Hopf algebra. Other variants are  possible as well, subject only
to the existence of natural examples.

Next we consider what should be the right concept of `unitary' representation
for a quasi-* Hopf algebra. The minimum definition, which is familiar for
groups and usual $*$-quantum groups, is a vector space $V$ on which the Hopf
algebra is represented, and a sesquilinear form $(\ ,\ )_V:V\tens V\to\C$
(antilinear in its first input) such that
\eqn{unitary}{(h^*\la v,v')_V=(v,h\la v')_V.}
for all $v,v'\in V$. We do not insist for the moment on conjugation-symmetry,
non-degeneracy and positivity of the sesquilinear form. While this definition
looks innocent enough, the new feature of quasi-* Hopf algebras (in contrast to
groups and usual $*$-quantum groups) is that such a definition is  not
naturally closed under tensor product. I.e. there appears to be no general way
to combine $(\ ,\ )_V$ and $(\ ,\ )_W$ to define a new sesquilinear form $(\ ,\
)_{V\tens W}$ such that the tensor product representation $V\tens W$ obeys the
same condition (\ref{unitary}). The problem does not show up for any one
unitary representation but only when we try to define consistently the category
of all unitaries. We formulate now a more correct notion which {\em is} closed
under tensor product and then explain why it cannot be restricted to unitaries
for fundamental reasons.

\begin{defin} A {\em mutually adjoint pair} of representations of a quasi-$*$
Hopf algebra $H$ is a vector space $V$, a sesquilinear form $(\ ,\ )_V:V\tens
V\to\C$ and {\em two} actions $\la,\bar\la$ of the Hopf algebra, such that
\[ (h^*\bar\la v,v')_V=(v,h\la v')_V\]
for all $v,v'\in V$ and $h\in H$. A morphism between mutually adjoint pairs is
a pair of intertwiners $\phi,\psi:V\to W$, one for the $\bar\la$
representations and one for the $\la$ representations, such that
$(\phi(v),\psi(v'))_W=(v,v')_V$. A {\em quasiunitary} representation is an
adjoint pair for which $\bar\la$ and $\la$ are isomorphic representations.
\end{defin}

We recover the previous notion (\ref{unitary}) of unitarity as the diagonal
case where $\bar\la=\la$. In fact, since we do not demand that the sesquilinear
form is `conjugate symmetric', a quasiunitary representation also leads to a
unitary one by absorbing any isomorphism between $\bar\la,\la$ into the
sesquilinear form. The more general setting of adjoint pairs has a natural
tensor product. We give two (equivalent) descriptions of it.

\begin{propos} Let $H$ be a quasi-* Hopf algebra (or more generally, a
cocycle-* Hopf algebra). Then two mutually adjoint representations $V,(\ ,\
)_V$ and $W, (\ ,\ )_W$ have a tensor product
\[ V\tens W,\quad  (v\tens w,v'\tens w')_{V\tens W}=(\CR\umo\bar\la
v,v')_V(\CR\umt\bar\la w,w')_W\]
for all $v,v'\in V$ and $w,w'\in W$, where the action $\bar\la$ of $H$ in the
first input extends to tensor products with the opposite coproduct. I.e. we
regard $\bar\la$ as an $H^{\rm cop}$-module and $\la$ as an $H$-module. The
category of adjoint pair representations is monoidal.
\end{propos}
\proof If $V,W$ with their sesquilinear forms are two mutually adjoint pairs in
the sense of Definition~4.2 then
\align{&&\equad (h^*\bar\la(v\tens w),v'\tens w')_{V\tens W}=(h^*\t\bar\la
v\tens h^*\o\bar\la w,v'\tens w')_{V\tens W}\\
&&=((\CR\umt h\o\CR\ut)^*\bar\la v\tens (\CR\umo h\t\CR\uo)^*\bar\la w,v'\tens
w')_{V\tens W}\\
&&=(h\o{}^*\CR\umt{}^*\bar\la v,v')_V(h\t{}^*\CR\umo{}^*\bar\la w,w')_W\\
&&=(\CR\umo\bar\la v,h\o\la v')_V(\CR\umt\bar\la w,h\t\la w')_W=(v\tens
w,h\la(v\tens w))_{V\tens W}}
as required, where we used the definition of $(\ ,\ )_{V\tens W}$, the action
of $H$ in its first input using the opposite coproduct, the reality assumption
on $\CR$ in Definition~2.6, the assumption that $V,W$ are adjoint pairs, the
usual action of $H$ in its second input and the definition of $(\ ,\ )_{V\tens
W}$ again. Hence the tensor product is also a mutually adjoint pair. We define
the tensor product of two morphisms to be their usual tensor product as linear
maps. This correctly connects the corresponding tensor product sesquilinear
forms because each morphism (as an intertwiner) commutes with the action of
$\CR^{-1}$. This makes the tensor product of adjoint pairs of representations a
functor from two copies to one copy of the category.

Moreover, this construction is associative by the usual vector space
associativity. Thus
\align{ &&\equad ((v\tens w)\tens z,(v'\tens w')\tens z')_{(V\tens W)\tens
Z}=(\CR\umo\bar\la(v\tens w),v'\tens w')_{V\tens W}(\CR\umt\bar\la z,z')_Z\\
&&=(\CR'\umo\CR\umo\t\bar\la v,v')_V(\CR'\umt\CR\umo\o\bar\la
w,w')_W(\CR\umt\bar\la z,z')_Z\\
&&=(\CR\umo\bar\la v,v')_V(\CR'\umo\CR\umt\t\bar\la
w,w')_W(\CR'\umt\CR\umt\o\bar\la z,z')_Z\\
&&=(\CR\umo\bar\la v,v')_V(\CR\umt\bar\la(w\tens z),w'\tens z')_{W\tens
Z}=(v\tens (w\tens z),v'\tens (w'\tens z'))_{V\tens(W\tens Z)}.}
This makes the category of adjoint pairs of representations monoidal in the
sense of \cite{Mac:cat}.
\endproof

So even if $\bar\la=\la$ for some representations, as we take tensor products
of them the composite  $\bar\la,\la$ will begin to diverge when our coproduct
is truly non-cocommutative. They need not even remain isomorphic. Clearly, one
does not see this interesting phenomenon for groups or enveloping algebras. For
a quasi-$*$ Hopf algebra the opposite coproduct is nevertheless twisting
equivalent to the conjugate coproduct $\bar\Delta$, and since twisting does not
change the category of representations up to equivalence we can work
equivalently with the first input of $(\ ,\ )_V$ living in the category of
representations of $\bar H$ instead. Hence we can `neutralise' the above
appearance of the cocycle in the tensor product of our sesquilinear forms
provided we put it into a more complicated form for the tensor product of
$\bar\la$.

\begin{propos} Let $H$ be a quasi-* Hopf algebra (or more generally, a
cocycle-* Hopf algebra). Then two adjoint pairs of representations $V,(\ ,\
)_V$ and $W, (\ ,\ )_W$ have tensor product $V\tens W$,
$(v\tens w,v'\tens w')_{V\tens W}=(v,v')_V(w,w')_W$  where  the action in the
first input extends to tensor products using the conjugate coproduct of $H$.
I.e. we regard $\bar\la$ as a $\bar H$-module and $\la$ as an $H$-module.
\end{propos}
\proof This is entirely equivalent to the preceding proposition, and the proof
is similar. The extension of $\bar\la$ is  $h\bar\la(v\tens w)=\CR\umo
h\t\CR\ut\bar\la v\tens \CR\umt h\o\CR\ut\la w$ in terms of the coproduct of
$H$. The tensor product representations $V\tens W$ are not the same as before
but equivalent via the morphism $\CR^{-1}\bar\la,\CR^{-1}\la$ mapping the
adjoint pair on $V\tens W$ as defined in Proposition~4.3 over to the adjoint
pair on $V\tens W$ as defined presently. This provides a monoidal functor
between the two categories. \endproof

Both of these tensor products of mutually adjoint pairs can be useful. To be
concrete we focus now on the latter formulation since it allows us to keep the
usual tensor product of our sesquilinear forms. It also makes clear that our
problem of $\bar\la$ and $\la$ diverging would not arise if $H$ were a Hopf
$*$-algebra in the usual sense.
We consider now how to construct adjoint pairs of representations from actions
on $*$-algebras. If they are also isomorphic then we will have a unitary
representation after adjusting the sesquilinear form. We begin with an
elementary lemma which is natural but not {\em quite} what we need for our
examples from bosonisation. After that, we will modify it to accommodate the
case of interest.

\begin{lemma} Let a $*$-algebra $C$ be acted upon by a quasi-$*$ Hopf algebra
$H$ by an action $\la$ making it an $H$-module algebra. Let $\bar\la$ be the
{\em conjugate} representation making $C$ a $\bar H$-module algebra such that
\[ (h\bar\la c)^*=(\bar S h)^*\la c^*,\quad (h\la c)^*=(Sh)^*\bar\la c^*\]
for all $h\in H$ and $b,c\in C$. If $\phi:C\to \C$ is a $\la$-invariant linear
functional then the sesquilinear form
\[ (b,c)_{\phi}=\phi(b^*c)\]
makes $\bar\la,\la$ a mutually adjoint pair in the sense of Definition~4.3.
Moreover, $\bar\phi=\overline{\phi((\ )^*)}$ is $\bar\la$ invariant and
$\overline{(b,c)_{\phi}}=(c,b)_{\bar\phi}$.
\end{lemma}
\proof Firstly, $\la$ and either one of the displayed conditions stated
determines $\bar\la$ (we give both forms to maintain the symmetry; they are
equivalent). Then
$hg\bar\la c=((S^{-1}(h^*))(S^{-1}(g^*))\la c^*)^*=(S^{-1}(h^*)\la(g\bar\la
c)^*)^*
=h\bar\la(g\bar\la c)$ so we have an action, and
\[ h\bar\la(bc)=(S^{-1}(h^*)\la(c^* b^*))^*=(S^{-1}(h^*\o)\la
b)^*(S^{-1}(h^*\t)\la c^*)^*=(h\baro\bar\la b)(h\bart\bar\la c)\]
so we have covariance with respect to $\bar\Delta h=h\baro\tens h\bart$. Also,
given $\phi$ it is clear that $H$-invariance means
\eqn{phiinv}{\phi((S^{-1}h\la b)c))=\phi(h\t\la((S^{-1}h\o\la
b)c))=\phi((h\t\la (S^{-1}h\o\la b))(h\th\la c))=\phi(b(h\la c))}
for all $h\in H$ and $b,c\in C$. We will later need to consider this equation
also with $\phi$ only partially invariant. So $(h^*\bar\la
b,c)=\phi((h^*\bar\la b)^*c) =\phi((S^{-1}h\la b^*)c)=\phi(b^*(h\la c))=(b,h\la
c)$ as required, using our definition of $\bar\la$ as conjugate to $\la$. The
last line of the lemma maintains the symmetry between $\bar\la,\la$ and follows
at once from their mutually conjugate relationship as stated. \endproof

In the setting of actions on $*$-algebras, our consideration of pairs of
actions  $\bar\la,\la$ is the same as considering either one (since one
determines the other). If, on the other hand, $\bar\la,\la$ are given to us
from some other source then the displayed condition in the  lemma becomes a
non-trivial constraint that the two actions are mutually conjugate. The lemma
generalises  standard considerations for unitary actions of Hopf $*$-algebras
to our setting of pairs of actions of quasi-$*$ Hopf algebras. If
$\phi=\bar\phi$ then the resulting sesquilinear form is conjugate symmetric,
which is again the usual case. One can then proceed to construct a Hilbert
space from this data in the standard way.

Unfortunately, our basic examples such as the fundamental and conjugate
fundamental representation of the  quasi-* Hopf algebras obtained by
bosonisation, do not necessarily fit into this standard setting and we have to
modify it. The variant we need is the following. It generally destroys,
however, the conjugation-symmetry of the sesquilinear form and hence forces us
to modify this axiom of Hilbert space theory.

\begin{lemma} Let $(C,\star)$ be a $*$-algebra acted upon by a quasi-$*$ Hopf
algebra $H$ by actions $\la,\bar\la$ as in Lemma~4.6 but now mutually conjugate
in the sense $(h\bar\la c)^\star=S(h^*)\la c^\star$ for all $c\in C$ and $h\in
H$. Let $\theta:C\to C$ be a linear automorphism such that $\theta(h\la
c)=S^{-2}h\la \theta(c)$ for all $h,c$. Then
$(b,c)^\theta_\phi=\phi((\theta(b^\star)c)$ makes $\bar\la,\la$  mutually
adjoint. The conjugate construction has just the same form with
\[ \bar\phi=\overline{\phi((\ )^\star)},\quad
\bar\theta=\star\circ\theta^{-1}\circ\star,\quad \overline{(c,b)^\theta_\phi}
=(b,c)^{\bar\theta}_{\bar\phi\circ\bar\theta^{-1}}.\]
\end{lemma}
\proof This is a variant of the preceding proposition. We have
$(h^*\bar\la b,c)=\phi(\theta((h^*\bar\la b)^\star)c)=\phi(\theta(Sh\la
b^\star)c)=\phi((S^{-1}h\la\theta(b^\star))c)=(b,h\la c)$ much as before. For
the conjugate we deduce $\bar\theta(h\bar\la c)=(\bar
S^{-2}h)\bar\la\bar\theta(c)$ when this is defined as stated. Hence the
conjugate construction has the same form. \endproof

If we had $\bar\theta=\theta$ then we would be able to redefine
$\theta\circ\star$ as a new $*$-structure and  precisely return to the setting
of Lemma~4.6. This is, however, not necessarily the case for the examples
coming from bosonisation of braided groups which interest us here. The origin
of the problem (already noted in \cite{Ma:star}) is that the duality pairing
(\ref{dual*}) for $*$-braided groups does not preserve the unitarity condition
(\ref{act*}) for the action of the background quantum group generating the
category. Indeed, the natural `unitarity' of the action on the $*$-braided
group $C$ dual to $B$ where the action obeys (\ref{act*}), is instead
\eqn{actstar}{(h\la c)^\star=S(h^*)\la c^\star}
for $c\in C$ and $h$ in the background quantum group. This is clear from
invariance of $\ev$ and computation of $\overline{\ev(h\la
b,c)}=\overline{\ev(b,(Sh)\la c)}$. When the duality pairing is degenerate then
$\star$ is not fully determined by
(\ref{dual*}) but we nevertheless keep (\ref{actstar}) as a reasonable
assumption compatible with the condition (\ref{act*}) which we assumed for the
action on $B$.

\begin{propos} Let $B\lbiprod H$ be the quasi-$*$ Hopf algebra constructed as
in Corollary~2.7, where $H$ is a real-quasitriangular Hopf algebra acting
unitarily on $*$-braided group $B$. Let $C$ a $*$-braided group dual to $B$ on
which $H$ acts as in (\ref{actstar}). The fundamental and conjugate fundamental
representation from Corollaries~2.2 and~2.4 of $B\lbiprod H$ on $C$ are then
mutually conjugate in the sense of Lemma~4.7:
\[ (x\bar\la c)^\star=S(x^*)\la c^\star,\quad (x\la c)^\star=\bar S(x^*)\bar\la
c^\star\]
for all $x\in B\lbiprod H,c\in C$. Here $(\ )^\star$ is the $*$-structure on
$C$ characterised by (\ref{dual*}) and $S,\bar S$ the antipodes for $B\lbiprod
H$.
\end{propos}
\proof We compute the first of these for $b\in B$ from the definition in
Corollary~2.4 as
\align{&&\equad (b\bar\la
c)^\star=c\Bt{}^\star\ev(b^*,c\Bo{}^\star)=\CR\ut\CR\umt\la c^\star\Bo\ev(\und
S^{-1}\und S(b^*),\CR\uo\CR\umo\la c^\star\Bt)\\
&&=(S\CR\ut)\CR\umt\la c^\star\Bo\ev(\und S^{-1}(\CR\uo\la \und
S(b^*)),\CR\umo\la c^\star\Bt)\\
&&=(S\CR\ut)\la((\CR\uo\la \und S(b^*)\la
c^\star)=\left(((S\CR\ut\t)\CR\uo\la\und S(b^*))S\CR\ut\o\right)\la c^\star\\
&&=(\cu\CR\uo\la\und S(b^*))S\CR\ut\la c=S(b^*)\la c^\star}
where the first equality is the definition of $\bar\la$, while the second
inserts $\CR\CR^{-1}$ and   uses one of the $*$-braided group axioms
(\ref{bg*}) to move $\star$ onto $c$. For the third equality we replace $\CR$
by $(S\tens S)(\CR)$ and move $S\CR\uo$ over to the left hand input of $\ev$
(by its $H$-invariance). It passes through $\und S^{-1}$ also, since this is
$H$-covariant. The fourth equality recognises the fundamental representation
from Corollary~2.2. The fifth uses that that this too is a representation of
$B\lbiprod H$, allowing us to move the $S\CR\ut$ to the left via its relations
(\ref{bos}). The sixth computes this via the coproduct action for $\CR$, which
gives the element $\cu$. We then recognise the antipode $S$ if $B\lbiprod H$
from (\ref{bos}). The proof for
$(b\la c)^\star$ is analogous, giving this time $\bar S$ from Proposition~2.3.
For the $x\in H$ we know that its actions $\bar\la$ and $\la$ coincide with the
action of $H$ on $C$ as a braided group, so this case is (\ref{actstar})
characterised by the braided group duality (\ref{dual*}) as explained above.
\endproof

Hence Lemma~4.7 is the variant of the general theory for quasi-$*$ Hopf
algebras that we need for our bosonisation examples $B\lbiprod H$ with $\la$
the fundamental representation. Next we consider the functional $\phi$. In the
present case it is clear from the definition of $\la$ in Corollary~2.2 that
when the duality pairing $\ev$ is non-degenerate, a $B\lbiprod H$-invariant
complex linear map $C\to \C$ means nothing other than an $H$-covariant and
right-invariant braided integration $\int{}$ on the braided group $C$.
Likewise, a $\bar\la$-invariant linear map means nothing other than an
$H$-covariant and left-invariant integration $\bar{\int}{}$. These maps are
characterised as morphisms $C\to\C$ by
\eqn{braint}{ (\int c\Bo)c\Bt=\int c,\quad c\Bo\int_L c\Bt=\int_L c}
for all $c\in C$, respectively. It is clear from (\ref{bg*}) that for
$*$-braided groups a morphism $\int{}$ is a right integral {\em iff}
$\bar{\int}=\overline{\int(\ )^\star}$ is a left integral. At least for
finite-dimensional and other reasonable braided groups, the left and right
integrals (if they exist) are unique up to normalisation. So we could take $
\int_L=\bar{\int}$ without significant loss of generality.

\begin{propos} In the setting of Proposition~4.8, let $\int{},\int_L{}$ be
right-i, left-invariant integrations for the braided group $\C$. Then the
conjugate/fundamental representations $\bar\la,\la$ are mutually adjoint with
respect to the sesquilinear form $(\ ,\ )$ and $\la,\bar\la$ with respect to
the sesquilinear form $(\ ,\ )^L$ defined by
\[ (b,c)=\int \theta_{\cv}(b^\star)c,\quad
(b,c)^L=\int_L\theta_{\cu}^{-1}(b^\star)c \]
for all $b,c\in C$. Moreover, $\overline{(c,b)}=(b,c)^L$ when
$\int_L=\bar{\int}\circ\theta_\cu$.
\end{propos}
\proof One can verify directly that the first sesquilinear form makes
$\bar\la,\la$ mutually adjoint. We have already done the work however, and
deduce it as follows. From the mutual conjugation property in Proposition~4.8
and the assertion in Corollary~2.4 that $\bar\la,\la$ are intertwined by $\und
S$, we easily deduce that
\eqn{S^2la}{\und S^2(x\la c)=((S\circ *)^{-2}x)\la \und S^2c,\quad \und
S^2(x\bar\la c)=((S\circ *)^{-2}x)\bar\la \und S^2c}
for all $x\in B\lbiprod H$ and $c\in C$. From this, (\ref{S*}), and
Proposition~4.2 it is clear that
\eqn{thetav}{\theta_\cv(x\la c)=(S^{-2}x)\la\theta_\cv(c),\quad
\theta_\cv(x\bar\la c)=(S^{-2}x)\bar\la\theta_\cv(c).}
Hence we have the required automorphism $\theta=\theta_\cv$ for Lemma~4.7. By
similar considerations to those above, we deduce equally well
\eqn{thetau}{\theta_\cu(x\la c)=(\bar S^2 x)\la\theta_\cu(c),\quad
\theta_\cu(x\bar\la c)=(\bar S^2 x)\bar\la\theta_\cu(c)}
for all $x\in B\lbiprod H$ where $\theta_\cu=\star\circ\theta_\cv\circ\star$.
The latter follows from $\cv^*=\cv$, $S\cv=\cu$ and (\ref{bg*}). This gives the
conjugate construction for the sesquilinear form $(\ ,\ )^L$ obeying $(h^*\la
b,c)^L=(b,h\bar\la c)^L$. Then
\align{&&\equad \overline{\int\theta_\cv(b^\star)c}=\bar{\int}
b^\star\theta_\cu(c)=\bar{\int}(\cv\la b^\star)\und
S^2c=\bar{\int}\theta_\cu(\theta_{\cu}^{-1}(b^\star)c)}
explicitly relates the two constructions as stated. The third expression means
that we can also write $\overline{(c,b)}=(\und S^{-2}b,\und S^2c)$ if we
replace $\int{}$ on the right by $\bar{\int}{}$. \endproof

We have $\bar{\int}\circ\theta_\cu=\bar{\int}$ as well in reasonable cases
where the left-invariant integration is unique up to scale. This achieves the
task of making our fundamental and conjugate fundamental representations from
Corollaries~2.2 and~2.4 mutually adjoint. Since we already know  that they are
equivalent with intertwining operator $\und S$, we can absorb this too into the
sesquilinear form. Then
\eqn{sesquni}{ (b,c)^U_{\phi}=\int(\cv\la\und Sb^\star)c}
for $b,c\in C$ defines a sesquilinear form with respect to which the
fundamental representation of $B\lbiprod H$ in Corollary~2.2 is unitary. In our
interpretation as Poincar\'e quantum group represented on spacetime, it
corresponds to building a parity operator into the $L^2$ inner product. With
respect to such a non-local and non-symmetric sesquilinear form we would have
$\del^i$ self-adjoint (i.e., symmetric) rather than the usual anti-selfadjoint.
The sesquilinear form in Proposition~4.9 by contrast becomes the usual $L^2$
inner product.

As explained above, it is natural in these constructions to consider
$\theta_\cv\circ\star$ as a possible second $*$-structure on $C$. It is an
antilinear anti-algebra homomorphism but
\eqn{theta2}{ (\theta_\cv\circ\star)^2=\cu\cv\la\und S^4}
is not necessarily the identity, though it is in some cases. In the present
setting, $\cu\cv$ is a central element of our quasitriangular Hopf algebra. One
could also build in its square root when this exists (the so-called ribbon
element\cite{ResTur:rib}) in order to reduce the contribution of $\cv$ to
$\theta^2$. Partly in this direction, it is easy to see that if (as in the
ribbon case) we have\cite{Dri:alm} a group-like element $\sigma$ implementing
$S^2$ and such that $\sigma^*=\sigma$ then any $*$-braided group $(C,\star)$ in
the category of $H$-modules has another $*$-structure $c^*=\sigma^{-1}\la
c^\star$. Moreover, this converts (\ref{actstar}) over to the more standard
unitarity condition  (\ref{act*}). If we use this second $*$-structure then the
sesquilinear form in  Proposition~4.8 becomes
\eqn{ribbonL2}{ (b,c)=\int \theta_\nu(b^*)c}
where $\nu=(\cv\cu)^{\h}$ is the ribbon element. This is a purely cosmetic
change.

We conclude this abstract section with an example of the above theory which is
somewhat different from the inhomogeneous quantum groups  of primary interest
in the present paper. Namely, we showed in \cite{Ma:skl} how to view Drinfeld's
quantum double $D(U_q(\cg))$ as a bosonisation $BU_q(\cg)\lbiprod U_q(\cg)$,
where $U_q(\cg)$ is  from \cite{Dri}\cite{Jim:dif} and
$BU_q(\cg)$ is its associated braided group from \cite{Ma:skl}. We consider the
standard Lie algebra deformations where there is an R-matrix form
$\vecl^\pm,\vect$ as in \cite{FRT:lie} of the quantum group and its dual $G_q$.
The dual of $BU_q(\cg)$ is the matrix braided group $BG_q$ obtained as a
corresponding quotient of the braided matrices $B(R)$\cite{Ma:exa}. We take
$q^*=q$ and the standard compact real form of the quantum groups. As a
$*$-algebra $BU_q(\cg)$ coincides with $U_q(\cg)$ and we take matrix generator
$\vecm=\vecm^+S\vecm^-$ forming a $*$-braided group. The braided coproduct is
$\und\Delta\vecm=\vecm\tens\vecm$ and the $*$-structure is the Hermitian one.
We refer to \cite{Ma:skl} and \cite{Ma:mec} for full details of this
bosonisation.

\begin{propos} The quantum double $D(U_q(\cg))$ in the bosonisation form
$BU_q(\cg)\lbiprod U_q(\cg)$ is a quasi-$*$ Hopf algebra by Corollary~2.7. The
fundamental and conjugate fundamental representations are
\[ \vecl^+_1\la\vecu_2=R_{21}^{-1}\vecu_2R_{21},\quad
\vecl^-\la\vecu_2=R^{-1}\vecu_2 R,\quad \vecm_1\la R
\vecu_2=R\vecu_2R_{21}R,\quad \vecm_1\bar\la\vecu_2=R_{21}R\vecu_2\]
on $BG_q$ with matrix generator $\vecu$. Moreover, $\theta_\cv\circ\star=\und
*$ is the standard Hermitian $*$-structure on $BG_q$ and
\[ (b,c)=\int b^{\und *}\und\cdot c=\int \Theta(b^*)c,\quad \Theta(c)=c\th
\CR(S^2c\t, (Sc\o)c\fo)\]
where $\und\cdot$ is the product in $BG_q$. The second expression computes this
further in terms of $b,c\in G_q$ and its usual coproduct, antipode, $*$,
dual-quasitriangular structure\cite{Ma:euc} and Haar measure\cite{Wor:com}
$\int$. Here we identified $BG_q=G_q$ as coalgebras by
transmutation\cite{Ma:eul}\cite{Ma:bg}.
\end{propos}
\proof The bosonisation and its $*$-structure for real $q$ is described in
detail in \cite{Ma:mec} (as well as for $q$ modulus 1, which we do not consider
here). The braided group duality pairing between $B=BU_q(\cg)$ and $C=BG_q$ is
also given there, as $\ev(b,c)=\<Sb,c\>$ where the right hand side views $b\in
H=U_q(\cg)$ and $c\in G_q$ since they coincide as linear spaces ($B$ has a
modified coproduct and $C$ a modified product, making them braided groups in
the category of $H$-modules by the quantum adjoint and coadjoint action
respectively). Then the fundamental representation from Corollary~2.2 is
\[ b\la c=\CR\umt\la c\o\ev(\und S^{-1}b,\CR\umo\la
c\t)=\<\CR\ut,Sc\o\>c\t\<\CR\uo\la b,c\th\>\]
by a computation similar to that for braided right-invariant vector fields in
\cite[Prop. 6.2]{Ma:lie}.
The conjugate fundamental from Corollary~2.4 is just $b\bar\la
c=\ev(b,c\o)c\t=\<Sb,c\o\>c\t$. Putting in the matrix coproduct for the
generator $\vect\in G_q$ and a standard computation from (\ref{lpm}) gives the
form for these actions as stated. The action of the  $U_q(\cg)$ part is the
same for $\la$ and $\bar\la$ and is easily computed\cite{Ma:exa} by evaluation
of $\vecl^\pm$ against the quantum adjoint coaction.

Next, the operation $\star$ on $BG_q$ dual to that on $BU_q(\cg)$ in the sense
of (\ref{dual*}) comes out as
\[ c^\star=S^{-3}(c^*),\quad i.e.,\quad  \theta_\cv(c^\star)=\sigma^{-1}\la
c^\star=(Sc)^*\]
in terms of the usual $*$ and antipode $S$ of $G_q$. Here $\la$ is the quantum
coadjoint actions hence it is immediate that  that $\sigma^{-1}\la=S^2$. With
rather more work, we may use the formula in \cite{Ma:bg} for the braided
antipode $\und S$ of $BG_q$ in terms of $G_q$ to obtain $\theta_\cv=S^2$ as
well. Since $S^2$ is an automorphism of the quantum group $G_q$  it necessarily
induces an  automorphism of the associated $BG_q$. In the same way,
$\theta_\nu=\id$. We recognise the combination $\und *=(S\ )^*$ (which we
underline to keep it distinct from the $*$ of $G_q$) as the Hermitian
$*$-braided group structure for $BG_q$ introduced in \cite{Ma:mec}. Hence
\align{(b,c)\equad &&=\int \theta_\cv(b^\star)\und\cdot c=\int b^{\und
*}\und\cdot c=\int ((Sb)^*)\t c\t \CR(Sc\o,(S((Sb)^*)\o)((Sb)^*)\th)\\
&&=\int (\CR\ut\la(Sb)^* )c\t \<\CR\uo,Sc\t\>}
where the third equality expresses the product $\und\cdot$ of $BG_q$ in terms
of $G_q$ using the transmutation formula in \cite{Ma:bg}. It is useful (but not
necessary) to make this conversion because the Haar weights for compact quantum
groups are already known\cite{Wor:com}, in some cases quite explicitly. Since
the $BG_q$ has the same coalgebra, we use the same $\int$. It is a morphism
$C\to\C$ because it is both left and right invariant on $G_q$. Invariance also
gives the form for $(\ ,\ )$ stated. This, and the actions $\la,\bar\la$ make
sense over $\C$ and do not require formal powerseries. \endproof

The bosonisation here for the simplest case $BU_q(su_2)\lbiprod U_q(su_2)$ was
studied in detail in \cite{Ma:mec} so we do not repeat this here. The new
feature is that we know now that its $*$-algebra, which  we developed there as
q-deformed Mackey quantisation of a particle moving on the mass hyperboloid in
$q$-Minkowski space, is a quasi-$*$ Hopf algebra. Moreover, the fundamental
representation $\vecm\la=\del$ extends as a `matrix braided-derivation' and was
computed (in a right hand version) for the $BSU_q(2)$ case in \cite{Ma:lie}, to
which we refer for details. These derivatives form braided left-invariant
vector fields on the braided group and realise the matrix braided-Lie algebra
$gl_{2,q}$\cite{Ma:lie}. From either point it view it is natural to consider
functional analysis on $BG_q$ as the algebra of `co-ordinate functions' of the
braided group. The above theory now provides the sesquilinear form via the
integral which to lowest order (carrying out required the transmutation
explicitly) comes out as
\ceqn{intBSU2}{\int 1=1,\quad \int a^2=\int ab=\int ac=\int b^2=\int bd=\int
c^2=\int cd=0,\quad\int ad={1\over 1+q^2},\quad\int d^2={1-q^{-2}\over 1+q^2}}
in terms of the generators $\pmatrix{a&b\cr c&d}$ and the product of
$BSU_q(2)$. The integrals of the other quadratic expressions are determined by
the relations of $BSU_q(2)$. Our left-invariant vector fields and their
conjugates are adjoint with respect to the corresponding sesquilinear form
$(b,c)=\int b^*c$ etc., where  $\pmatrix{a^*&b^*\cr c^* & d^*}=\pmatrix{a& c\cr
b& d}$ is the Hermitian $*$-braided group structure of $BSU_q(2)$. This
demonstrates the possibility of a braided approach to $q$-harmonic analysis on
braided versions of compact quantum groups, to be considered elsewhere. Note
that in this family of examples the sesquilinear form $(\ ,\ )$ clearly remains
conjugate-symmetric, i.e. we are in the usual setting as in Lemma~4.6 after
redefining $*$. We do lose, however, positive definiteness.

\section{Differential representation on spacetime and concluding remarks}

In this section we will consider our above results apply to the fundamental and
conjugate fundamental representations of the inhomogeneous quantum groups from
Section~3. We are now in a position to understand the subtle role of the $*$
structure, which is obviously crucial for the interpretation of these quantum
groups as, for example, the $q$-Poincar\'e group in $q$-deformed geometry.

We begin with the fundamental and conjugate fundamental representations
themselves on the `spacetime' braided group $\Vhaj(R',R)$ with co-ordinates
$\{x_i\}$ dual to the linear `momentum' part of the inhomogeneous quantum
group. The existence of a general fundamental representation on `spacetime' was
the main result of the braided approach in \cite{Ma:poi} in the form of a
`rotation+ translation' coaction of the dual quantum groups. Evaluating against
this coaction (or from Corollary~2.2),  we obtain at once
\eqn{poinact}{\vecl^+_1\la \vecx_2=\vecx_2
\lambda R_{21},\quad \vecl^-_1\la\vecx_2=\vecx_2 \lambda^{-1}R^{-1},\quad
\lambda^\xi\la x_i=\lambda x_i,\quad p^i\la x_j=-\delta^i{}_j}
which we have used already in (\ref{lpmcov}),(\ref{lmcov}) in Section~3. This
then extend to products as a module algebra.  Explicitly:

\begin{propos} The inhomogeneous quantum group $V(R',R)\lbiprod \widetilde{H}$
in the setting of Section~3 acts covariantly on the braided covector algebra
$\Vhaj(R',R)$ with generators $\vecx=\{x_i\}$ by the fundamental and conjugate
fundamental representations in Corollaries~2.2 and~2.4 as
\cmath{
\lambda^\xi\la(x_{i_1}\cdots x_{i_m})=\lambda^m x_{i_1}\cdots x_{i_m} \\
l^+{}^i{}_j\la (x_{i_1}\cdots x_{i_m})=\lambda^m x_{j_1}\cdots
x_{j_m}[1,m+1;R]^{ij_1\cdots j_m}_{i_1\cdots i_mj}\\
l^-{}^i{}_j\la (x_{i_1}\cdots x_{i_m})=\lambda^{-m} x_{j_1}\cdots
x_{j_m}[1,m+1;R_{21}^{-1}]^{ij_1\cdots j_m}_{i_1\cdots i_mj}\\
-p^i\la (x_{i_1}\cdots x_{i_m})=x_{j_2}\cdots x_{j_m}\left[m;
R^{-1}_{21}\right]^{ij_2\cdots j_m}_{i_1i_2\cdots i_m}\\
p^i\bar\la (x_{i_1}\cdots x_{i_m})=x_{j_2}\cdots x_{j_m}\left[m;
R\right]^{ij_2\cdots j_m}_{i_1i_2\cdots i_m}}
where $[1,m+1;R]=(PR)_{12}(PR)_{23}\cdots(PR)_{mm+1}$, and
$[m;R]=[1,2;R]+[1,3;R]+\cdots +[1,m;R]$ is the braided integer matrix in
\cite{Ma:fre}. The action of the $\lambda^\xi,\vecl^\pm$ generators is the same
for the two representations.
\end{propos}
\proof Here $\Vhaj(R',R)$ is the algebra $\vecx_1\vecx_2=\vecx_2\vecx_1R'$
forming a braided group \cite{Ma:fre}. Its braided group duality with the
braided vector algebra used for the linear part of the inhomogeneous quantum
group is $\ev(p^i,x_j)=\delta^i{}_j$. Since we know (by the bosonisation
theory) that the fundamental representation makes $\Vhaj(R',R)$ a module
algebra, the action on products is then determined. The action of $\vecl^\pm$
on products is immediate from their matrix coproduct. The braided integer
matrix is a sum of the corresponding matrices, so it is clear by induction that
$p^i$ acts by such matrices given the action of $\vecl^\pm$ and the form of the
coproduct and conjugate coproduct. In the case of the conjugate coproduct the
action of $p^i$ is necessarily braided differentiation since the map in
Corollary~2.4 is exactly its definition in \cite{Ma:fre}.  Note that we  have -
sign in the first action $\la$ of $p^i$ from the braided antipode in
Corollary~2.2, but not in the conjugate action $\bar\la$ in Corollary~2.4.
\endproof

In particular, the covariant action of $p^i$ using the conjugate coproduct
$\bar\Delta$ is precisely the braided differentiation $p^i=\del^i$ as
introduced for general braided linear spaces in \cite{Ma:fre}. In fact, we know
this without computation because evaluation against the braided coproduct as in
the abstract definition of $\bar\la$ in Corollary~2.4 is precisely definition
of $\del^i$ in \cite{Ma:fre} as an `infinitesimal coaddition' in the braided
approach. The covariant action of $p^i$ using the original coproduct $\Delta$
is by `conjugate' derivatives $-p^i=\bar\del^i$ in which the role of $R$ is
replaced by $R_{21}^{-1}$ when extending to products. It is an infinitesimal
translation from the right and corresponds to the `right derivatives'
$\overleftarrow{\del}$ in \cite{Ma:star}, converted over to left-acting
derivatives by means of the braided antipode. This is the reason for the extra
$-$ sign in the action of $p^i$. The $\del$ obey a braided-Leibniz rule with
$\Psi^{-1}$ as explained in \cite{Ma:fre}, while the $\bar\del$ obey a braided
Leibniz rule with $\Psi$ as we have seen already in the proof of Corollary~2.2.
The reversed matching of $\bar\la$ with $\del$ etc., is an historical accident
reflecting the fact that the unbarred $\Delta$ and unbarred $\del$ are each
natural in their own settings.

One can also consider $x_i$ as an operator on $\Vhaj(R',R)$ by left
multiplication then the corresponding braided-Leibniz rules are expressed as
the commutation relations
\eqn{leib}{  \del_1\vecx_2-\vecx_2R_{21}\del_1=\id,\quad
\bar\del_1\vecx_2-\vecx_2R^{-1}\bar\del_1=\id.}
cf. specific examples in \cite{PusWor:twi}\cite{Kem:sym}\cite{OSWZ:def}, etc.
If we assume a quantum metric and lower indices by $\del_i=\eta_{ia}\del^a$ and
$\bar\del_i=\eta_{ia}\bar\del^a$ then these become
 \eqn{lowleib}{ \del_1\vecx_2-\lambda^{-2}\vecx_2\del_1R^{-1}_{21}=\eta,\quad
\bar\del_1\vecx_2-\lambda^2\vecx_2\bar\del_1R=\eta}
using the quantum metric identities (\ref{metricid}). These represent the
lower-index momentum generators $p_i$ for the two actions.  This is how the
constructive braided approach to differentials\cite{Ma:fre} recovers previous
approaches\cite{PusWor:twi}\cite{WesZum:cov} where examples of such commutation
relations were postulated as an ansatz or deduced from postulated relations
between differential forms within an axiomatic approach for these.

We obviously have similar formulae for the representations of the `spinorial'
forms $p^{i_0}{}_{i_1}$ in Sections~3.1 and~3.2. This is just a change of
notation from the vector form to the matrix form, $\del_I=\del^{i_0}{}_{i_1}$
and $\bar\del_I=\bar\del^{i_0}{}_{i_1}$, say. Clearly
\eqn{eucdif}{R\del_2\vecx_1R-\lambda^{-2}\vecx_1\del_2=R\eta_{21} R, \quad
\bar\del_1\vecx_2-\lambda^2R\vecx_2\bar\del_1R=\eta}
in the $\bar A(R)$ case where $\del,\bar\del$ obey the $\bar A(R)$ relations
$R_{21}\del_1\del_2=\del_2\del_1R$, etc., just because this is how ${\bf
R}',{\bf R}$ appear in this notation (\ref{eucmat}). Likewise, we have
\eqn{minkdif}{\del_2R_{21}\vecu_1R-\lambda^{-2}R^{-1}\vecu_1R\del_2=\eta\ut
R_{21}\eta\uo R,\quad R^{-1}\bar\del_1
R\vecu_2-\lambda^2\vecu_2R_{21}\bar\del_1R=R^{-1}\eta\uo
R\eta\ut}
in the $B(R)$ case where $\del,\bar\del$ obey the $B(R)$ relations
$R_{21}\del_1 R\del_2=\del_2 R_{21}\del_1
R$ etc., because this is how ${\bf R'},{\bf R}$ appear in this notation
(\ref{bramat}). Here $\eta=\eta\uo\tens\eta\ut$ is
$\eta_{IJ}=\eta^{i_0}{}_{i_1}{}^{j_0}{}_{j_1}$ as an element of $M_n\tens M_n$,
and the right hand sides are typically multiplies of it as well. We include
these formulae for completeness only; they are just the standard construction
(\ref{lowleib}) applied to the particular ${\bf R}',{\bf R}$ introduced in
\cite{Ma:exa}\cite{Mey:new}\cite{Ma:euc}.

Because the braided approach derives such relations from the braided coproduct
rather than imposing them axiomatically, we are now in a position to say more
about the derivatives $\del,\bar\del$ than is evident from the relations alone.

\begin{corol} The braided antipode $\und S(\vecx)=-\vecx$ of the braided group
$\Vhaj(R',R)$ intertwines the actions of $\del,-\bar\del$,
\[ \und S\del^i=-\bar\del^i\und S.\]
\end{corol}
\proof The diagrammatic proof is given in the proof of Corollary~4.6 and is a
result entirely as operators on the braided group $\Vhaj(R',R)$ in the setting
of \cite{Ma:fre}. For a direct proof the key fact is that $\und S$ extends to
products as a braided-anti-algebra homomorphism (\ref{braant}), giving it an
expression in terms of R-matrices\cite{Ma:poi} which intertwines the braided
integer matrices in Proposition~4.1. The operator $\und S$ is also covariant
under the background quantum group $H$\cite{Ma:poi} and hence intertwines its
action also. It also preserves degree, hence the action of $\lambda^\xi$.
\endproof

This general result contrasts with other approaches,
e.g.\cite{SchWes:def}\cite{OgiZum:rea} where it is sometimes possible to write
the $\bar\del$ as some non-linear function of the $\del$. We have not taken
this line here: in the braided approach the above corollary is more natural
since the braided antipode $\und S$ is inversion on the additive braided group,
i.e. the `braided parity operator' and plays such a role in important role in
numerous other constructions as well, such the braided Fourier theory in
\cite{LyuMa:bra}\cite{KemMa:alg}. In usual undeformed constructions we see it
merely as a minus sign, but in the braided case it extends as a non-trivial
operator.

Next, we consider $*$-structures. We have seen in Section~4 the need to
consider a different $*$-structure $\star$ on the `spacetime' braided group $C$
determined by duality (\ref{dual*}) with the linear `momentum' part $B$ of our
inhomogeneous quantum group. We have used for the latter the standard
$*$-structure which is characterised by the unitarity condition (\ref{act*})
with respect to the action of the background quantum group. Hence for the
spacetime $*$-braided group we need one which is characterised by
(\ref{actstar}), a condition which is different when $S^2\ne \id$. This
splitting into two $*$-structures (even when the braided groups are isomorphic
via the quantum metric) is therefore a new feature of $q$-deformation.  The
duality pairing of $*$-braided groups in the present linear case was studied in
\cite{Ma:star} from where we may deduce the required $\star$ appropriate to the
$*$ used for the momentum generators $p^i$ in Section~3. For the type I case
(as in Propositions~3.1 and~3.4) and the real type II case (used in
Proposition~3.7) they are
\eqn{newstar}{ p^i{}^*=\cases{\eta_{ia}p^a&real\ type\ I\cr  \eta^{\bar i
a}\eta_{ab}p^b&real\ type\ II},\quad x_i{}^\star=\cases{x_a\eta^{ia}&real\
type\ I\cr  x_b\eta_{\bar ia}\eta^{ab}&real\ type\ II}}
where $\eta$ is the quantum metric and in the type II case we also have an
involution $\bar{\ }$ on its indices. This corresponds via the quantum metric
to $p_i{}^*=p_{\bar i}$ as used in Proposition~3.7 in the form
$p^{i_0}{}_{i_1}$ Hermitian. With  $\star$ defined correctly, we see from
Proposition~4.8 that it connects $\del,\bar\del$. In the real type I setting of
Proposition~3.2, this is
\eqn{delstar}{ (\del^i f(\vecx))^\star=\lambda^{\xi}
S^2l^-{}^a{}_i\la\bar\del_a f(\vecx)^\star}
for all multinomials $f(\vecx)$. We used the antipode from Proposition~3.1. A
more intrinsic formulation (purely in terms of the linear braided group) is
possible\cite{Ma:star} if we use right-acting derivatives
$\overleftarrow{\del}$ in place of $\bar\del$. There are corresponding formulae
for the spinorial versions in Sections~3.1, 3.2 as well.

We also needed in Section~4 the braided group automorphisms
$\theta_\cv,\theta_\cu,\theta_\nu$. Evaluating the corresponding elements of
our background quantum group $\widetilde{H}$ etc., against the quantum matrix
transformation of the $x_i$ gives at once
\eqn{vectheta}{\theta_\cv(x_i)=x_a\cv^a{}_i,\quad
\theta_\cu(x_i)=x_a\cu^a{}_i,\quad \theta_\nu(x_i)=\lambda_\nu x_a ;\quad
\cv^i{}_j=\widetilde{R}^i{}_a{}^a{}_j,\quad
\cu^i{}_j=\widetilde{R}^a{}_j{}^i{}_a,\quad (\cu\cv)^\h=\lambda_\nu\id}
where $\widetilde{R}=((R^{t_2})^{-1})^{t_2}$ is the second inverse (here $t_2$
denotes transposition in the second matrix factor of $M_n\tens M_n$.) The
action of the dilaton with its contribution $\lambda^{\xi^2}$ to $\cv,\cu$ and
$\nu$ cancels the quantum group normalisation factor otherwise appearing in
these formulae. The $\theta_\nu$ applies in the case that the background
quantum group is ribbon and the matrix representation is irreducible, which is
typical in examples. The extension of the $\theta_\cv,\theta_\cu,\theta_\nu$ is
as algebra homomorphisms.

\begin{propos} For the examples $\bar M_q(2)\lbiprod \widetilde{U_q(su_2)\tens
U_q(su_2)}$ and $BM_q(2)\lbiprod \widetilde{U_q(su_2)\codcross U_q(su_2)}$ of
4-dimensional $q$-Euclidean and $q$-Minkowski Poincar\'e groups, we have
\[ \cv=\pmatrix{q^{-4}&0&0&0\cr 0&q^{-6}&0&0\cr 0&0&q^{-2}&0\cr
0&0&0&q^{-4}},\quad \cu=\pmatrix{q^{-4}&0&0&0\cr 0&q^{-2}&0&0\cr
0&0&q^{-6}&0\cr 0&0&0&q^{-4}}, \quad \lambda_\nu=q^{-4}.\]
Moreover, on the $q$-Euclidean and $q$-Minkowski spacetime co-ordinates we have
\[\theta_{\cv}((\vecx_1\cdots\vecx_n)^\star)
=(\vecx_1\cdots\vecx_n)^*\lambda_\nu^n.\]
\end{propos}
\proof The first part is best computed from (\ref{vectheta}) using the
appropriate ${\bf R}$ from \cite{Ma:euc}\cite{Mey:new}, but can also be done
directly from the $\cv,\cu,\nu$ elements in each copy of $U_q(su_2)$.
Moreover, the $\star$ structure on the $q$-spacetime co-ordinates in these two
cases
comes out from (\ref{newstar}) as
\eqn{eucminkstar}{\pmatrix{a^\star&b^\star\cr c^\star & d^\star}=\pmatrix{d& -q
c\cr -q^{-1}b& a},\quad \pmatrix{a^\star&b^\star\cr c^\star &
d^\star}=\pmatrix{a& q^2 c\cr q^{-2}b& d}}
for the algebras $\bar M_q(2)$ and $BM_q(2)$ respectively. From this we see at
once that $\theta_\cv(x_i^\star)=x_i^*\lambda_\nu$ where $*$ is the standard
$*$-structure on the spacetime co-ordinates (obeying the unitarity condition
(\ref{act*})) for the two cases. We denote the spacetime co-ordinates in both
cases by $x_i$ and the general form of $*$ is\cite{Ma:star}
\[ x_i^*\cases{x_a\eta^{ai}& Real\ type\ I\cr x_{\bar i}&Real\ type\ II}.\]
Our background quantum groups in these examples are ribbon  and
$c^\star=\sigma\la c^*$ where  $\sigma=\cu\nu^{-1}$ is Drinfeld's group-like
element\cite{Dri:alm} computed from that of $U_q(su_2)$. So we can also take
the point of view leading to (\ref{ribbonL2}) for these examples.  \endproof

We see that $\theta_\cv\circ\star$ is not an involution, but it is very close
to one, differing only by multiples in each degree. There is a similar form for
the $SO_q(n)$-covariant $\R_q^n$ spaces in \cite{FRT:lie} regarded as linear
braided groups. Finally, we need for our constructions an invariant integral.
There are two problems here, both of which can be addressed. The first is that
we cannot expect polynomials in the co-ordinate generators $x_i$ of on our
$q$-deformed linear spaces to be integrable; this can handled by defining
directly a Gaussian-weighted integral instead\cite{KemMa:alg}. The second
problem is that we cannot expect the integral to be invariant under the dilaton
part $\lambda^\xi$ of the inhomogeneous quantum group. This second problem is
dealt with by slightly generalising our Lemma~4.7 as follows: suppose that
$\phi$ is a linear functional $C\to \C$ which is invariant under some
sub-algebra of the quasi-$*$ Hopf algebra, and that the coproduct $\Delta
h=h\o\tens h\t$ of a given element of the latter can be written with the $h\t$
parts lying in the subalgebra. Then we can still conclude $(h^*\bar\la
b,c)_\phi^\theta=(b,h\la c)^\theta_\phi$ as before. This is evident from
(\ref{phiinv}) where we wrote the required step needed in Lemmas~4.6 and~4.7
quite explicitly. So we consider now a map $\int$ on $\Vhaj(R',R)$ which is
invariant under translation with respect to $\bar\del$ and under the background
quantum group without the dilaton. Hence it is invariant under the subalgebra
of the inhomogeneous quantum group without the dilaton. From the form of the
coproduct in Proposition~3.1 we see that we have the right form for all the
generators except the dilaton $\lambda^\xi$. So the fundamental and conjugate
fundamental representations are indeed mutually adjoint as regards the actions
of $\vecp,\vecl^\pm$  (but not of $\lambda^\xi$) with respect to
\[ (b,c)=\int (\lambda_\nu^\xi\la b^*)c.\]
Here the action of $\lambda_\nu^\xi$ is multiplication of $x_i$ by
$\lambda_\nu$ in the setting of Proposition~5.3 and similar examples.

More precisely, the appropriately covariant integration functional on general
braided linear spaces has been constructed in \cite{KemMa:alg} in a
Gaussian-weighted form. We need to give a right-invariant version appropriate
to $\bar\del$ rather than $\del$ as given there. Briefly, suppose formally that
 there is a Gaussian $\bar g$ solving the equation $\bar\del^i\bar g=-x_a
\eta^{ai} \bar g$ as a power-series. With the restrictions on $R,\eta$ in
\cite[Sec. V.1]{KemMa:alg} it takes the form of a $q$-exponential $\bar
g=e_{\lambda^2}^{-(1+q^2)^{-1}\vecx\cdot\vecx}$ which is central and invariant
under $\vecl^\pm$ and $*$. We do not need its precise form explicitly, however.
Instead, we define directly a linear functional $\CZ:C\to \C$ which plays the
role of the ratio $\CZ(f(\vecx))=\int f(\vecx)\bar g/\int \bar g$. We regard
the left hand side as a definition of the right hand side. The former, in turn,
is defined directly in terms of the R-matrix $R$ by means of induction as
cf.\cite{KemMa:alg}
\ceqn{cz}{ \CZ[1]=1,\quad \CZ[x_i]=0,\quad \CZ[x_ix_j]=\lambda^{-2} \eta_{ab}
R^{-1}{}^a{}_j{}^b{}_i\\
\CZ[x_{i_1}\cdots x_{i_m}]=\sum_{r=0}^{m-2}\CZ[x_{i_1}\cdots
x_{i_r}x_{a_{r+3}}\cdots x_{a_m}]\CZ[x_{i_{r+1}}x_{a_{r+2}}]
[r+2,m;R_{21}^{-1}]^{a_{r+2}\cdots a_m}_{i_{r+2}\cdots
i_m}\lambda^{-2(m-2-r)}.}
We refer to \cite{KemMa:alg} for the detailed derivation (for the left-handed
case appropriate to  $\del$). For the example of $SO_q(n)$-covariant quantum
planes such Gaussian-weighted integrations are known by other more explicit
calculations as well\cite{Fio:sym}. We conclude in particular that
$\del,\bar\del$ are mutually adjoint with respect to the sesqulinear form
defined implicitly by
\eqn{braL}{{(b(\vecx),c(\vecx)\bar g)\over(1,\bar g)}={\int
\lambda_\nu^{|b|}b^*c \bar g\over \int \bar g}=\lambda_\nu^{|b|}\CZ(b^*c)}
where $b,c$ are multinomials in $x_i$ and $|\ |$ is the degree. More precisely,
we take the right hand side directly as a definition of a Gaussian-weighted
sesquilinear form $\CZ(b,c)=\lambda_\nu^{|b|}\CZ(b^*c)$ even when $\int{}$ and
$\bar g$ are not defined.
The adjointness property then becomes
\ceqn{adjdel}{\CZ((\vecl^\pm)^*\la b, c)=\CZ(b,\vecl^\pm\la c),\quad
\CZ((\del^i)^*b,c)=-\CZ(b,\bar\del^ic)+\CZ(b,\cdot\Psi(x_a\tens
c))\eta^{ai}\lambda^{2|c|}}
where $\cdot\Psi(x_{i_1}\tens x_{i_2}\dots x_{i_m})=x_{a_1}\cdots
x_{a_m}[1,m;R_{21}]^{a_m\cdots a_1}_{i_m\cdots i_1}$ from the standard braiding
in the covector algebra\cite{Ma:fre}, and $\del^i{}^*=\eta_{ia}\del^a$, etc. is
as for $p^i{}^*$ in (\ref{newstar}), depending on the case. The two terms on
the right come immediately from the computation of $\bar\del(c\bar g)$ using
the braided-Leibniz rule followed by the equation for the Gaussian.
Equivalently, we use the Poincar\'e algebra coproduct (\ref{vecpoienv}) and the
action of $\lambda^\xi,\vecl^-$ in Proposition~5.1. Finally, once
(\ref{adjdel}) is obtained it may be verified directly for our $q$-Minkowski
and $q$-Euclidean examples (at least to low order) on multinomials $b,c$ even
when $\int,\bar g$ are not given. A formal inductive proof is rather long and
will be considered elsewhere. Non-degeneracy is also clear for our standard
examples (and for generic $q$) since it holds for $q=1$. Finally, our standard
examples are unimodular which, combined with \cite[Sec. 4]{Ma:star}, tells us
that $\CZ=\overline{\CZ((\ )^*)}$ as one may verify directly for low order.
This means that our sesquilinear form is conjugation-symmetric in the deformed
sense
\eqn{conjsymZ}{
\overline{\CZ(c,b)}={\lambda_\nu^{|c|}\over\lambda_\nu^{|b|}}\CZ(b,c).}
One can also leave out $\lambda_\nu$ in the definition of $\CZ(\ ,\ )$ and have
the more standard conjugation-symmetry, but at the price of a spurious factor
$\lambda_\nu^{-1}$ on the left hand side of the adjointness property of $\del$
in (\ref{adjdel}). This adjointness, combined with Corollary~5.2 is the sense
in which the differential representation of the $q$-Poincar\'e algebra or
inhomogeneous quantum group is `unitary' in our braided approach. The
sesquilinear form $\CZ(\ ,\ )$ appears to be the appropriate starting point for
$q$-quantum mechanics and $q$-scaler field theory in this approach. It remains
to develop suitable tools (such as a braided version of Wick's theorem) for the
computation of these Gaussian-weighted integrals and `braided $L^2$ inner
products' in a more closed form. It also remains to consider the appropriate
formulation of completions, domains of operators etc. for the corresponding
functional analysis. This is a direction for further work.

\baselineskip 18pt
%\bibliographystyle{unsrt}
%\bibliography{biblio}

 \end{document}